\documentclass{amsart}

\usepackage{amssymb}
\usepackage{euscript}
\usepackage{hyperref}

\newcommand{\E}{E}
\newcommand{\F}{\mathcal{F}}           % Stable family of diagonal seminorm functions
\renewcommand{\H}{\mathcal{H}}
\newcommand{\I}{\EuScript{I}} % A stable ideal
\newcommand{\J}{\EuScript{J}} % Another stable ideal
\newcommand{\K}{\mathcal{K}} % The compact operators in AlgN
\renewcommand{\L}{\mathcal{L}} % Some other nest
\newcommand{\M}{\mathcal{M}} % Some other nest
\newcommand{\N}{\mathcal{N}} % The fixed continuous nest
\renewcommand{\O}{\Omega}      % A stable net
\renewcommand{\P}{\mathcal{P}} % The set of all families of intervals
\newcommand{\Q}{\Psi}        % Another stable net
\newcommand{\R}{\mathcal{R}}

\newcommand{\alg}{\operatorname{Alg}}
\newcommand{\rank}{\operatorname{rank}}
\newcommand{\dist}[2]{\operatorname{dist}(#1, #2)}
\newcommand{\bh}{B(\H)}              % B(H), the bounded operators on a hilbert space, seee \h
\newcommand{\norm}[1]{\left\|#1\right\|}        % The operator norm
\newcommand{\enorm}[1]{\|#1\|_\ess}  % The essential norm
\newcommand{\dsf}[2]{\left\|#2\right\|_{#1}} % A general diagonal seminorm function
\newcommand{\normn}[1]{\dsf{N}{#1}}   % A diagonal seminorm function at N
\newcommand{\emptydsf}{\dsf{N}{\;\cdot\;}}
\newcommand{\innerprod}[2]{\langle #1, #2 \rangle}   % inner product

\newcommand{\ess}{\mathrm{ess}}
\newcommand{\ad}[1]{\operatorname{Ad}_{#1}}

\newcommand{\kplus}{\K^+}   % One of the basic ideals of compact character
\newcommand{\kminus}{\K^-}  % One of the basic ideals of compact character
\newcommand{\eo}{\E_0}      % A standard example of a stable ideal
\newcommand{\ei}{\E_I}      % A standard example of a stable ideal
\newcommand{\st}{\;:\;}     % 'such that' separator in set definitions
\newcommand{\e}{\epsilon}   % A small epsilon
\newcommand{\wlim}{\operatorname{w-lim}}
\newcommand{\ceil}[1]{\left\lceil#1\right\rceil}  % Rounding up.
\newcommand{\floor}[1]{\left\lfloor#1\right\rfloor}  % Rounding down.
\newcommand{\cpt}[2]{#1^{(#2)}}
\newcommand{\integers}{\mathbb{Z}}
\newcommand{\naturals}{\mathbb{N}}
\newcommand{\reals}{\mathbb{R}}
\newcommand{\rationals}{\mathbb{Q}}

\theoremstyle{plain}
\newtheorem{thm}{Theorem}[section]
\newtheorem{prp}[thm]{Proposition}
\newtheorem{cor}[thm]{Corollary}
\newtheorem{lem}[thm]{Lemma}

\theoremstyle{definition}
\newtheorem{defn}[thm]{Definition}
\newtheorem{eg}[thm]{Example}

\theoremstyle{remark}
\newtheorem{rem}[thm]{Remark}

\begin{document}

\title{The Stable Ideals of a Continuous Nest Algebra II}
\author{John Lindsay Orr}
\date{\today}
\address{
  John L. Orr,
  Department of Mathematics,
  University of Nebraska -- Lincoln, 
  Lincoln, NE 68588-0130, 
  USA
}
\email{jorr@math.unl.edu}

\begin{abstract}
  We continue the study of the rich family of norm-closed, automorphism
  invariant ideals of a continuous nest algebra. First we present a unified
  framework which captures all stable ideals as the kernels of limits
  of diagonal compressions. We then characterize when two such limits give
  rise to the same ideal, and we obtain detailed information of the
  structure of sums and intersections of ideals.
\end{abstract}

\maketitle

\section{Introduction}

In \cite{Orr:StIdNeAl} we studied norm-closed, automorphism-invariant
ideals of continuous nest algebras, which we termed stable ideals.
Our study was motivated by the fact that a number of natural examples of stable
ideals of nest algebras have been identified and studied, such as the
compact operators of the algebra; the Jacobson radical of the algebra,
described in \cite{Ringrose:OnSoAlOp};
the strong radical (that is, the intersection of all maximal two-sided
ideals) which was described in \cite{Orr:MaIdNeAl}; and the ideals
studied by Erdos in \cite{Erdos:OnSoIdNeAl}. In addition,
other interesting ideals narrowly fail to be automorphism invariant
for subtle reasons, specifically Larson's ideal $\R^\infty$
\cite{Larson:NeAlSiTr}. In \cite{Orr:StIdNeAl} we were able to
find a complete description of all stable ideals of a continuous
nest algebra. In the present paper we continue this study. 

In the original description of stable ideals, we identified
a family of eleven exceptional minimal ideals, which we dubbed
ideals of compact character, and then gave a characterization
of all remaining stable ideals
(see Theorem~\ref{thm:orig-stable-ideal-theorem}, below).
In the present paper we start, in Theorem~\ref{thm:stable_ideals},
by giving a single, unified description
of the stable ideals, which brings together both the compact character
and non-compact character cases. This description generalizes other
characterizations given elsewhere of ideals of nest algebras as the kernels of
the limits of certain diagonal expectations.

As a result of Theorem~\ref{thm:stable_ideals}, we associate each
stable ideal with a net of generalized partitions of the identity.
In Section~\ref{sec:cofinal-nets} we investigate when two such
nets give rise to the same stable ideal, and find that the simplest
condition, that the two nets be cofinal in each other, is both necessary
and sufficient. In Section~\ref{sec:quotient-norms} we build on these
results to give natural formulas for the quotient norm by a stable
ideal, and in Section~\ref{sec:sums-of-ideals} we characterize the sums of stable ideals. The quotient
norm formulas substantially generalize the quotient norm formula
obtained by Ringrose for the Jacobson radical \cite{Ringrose:OnSoAlOp}
and similar formulas found by Erdos \cite{Erdos:OnSoIdNeAl} for related ideals.

It is only possible to get such detailed information about stable ideals
of continuous nest algebras because of the very clear understanding
we have of the automorphism groups of nest algebras in general. 
Two key results underlie all of the work in this paper:
First Ringrose~\cite{Ringrose:OnSoAlOp2} showed that isomorphisms
between nest algebras are necessarily spatial; that is, every
isomorphism $\Phi:\alg\L\rightarrow\alg\M$ can be expressed
as $\Phi=\ad{S}$, where $S$ is a bounded, invertible operator.
Then Davidson~\cite{Davidson:SiCoPeNeAl} gave a necessary and sufficient condition
for two nest algebras to be isomorphic in terms of order-dimension
invariants of the nests:

\begin{thm}\label{thm:similarity}
  Let $\L$ and $\M$ be two nests on separable Hilbert spaces
  and suppose there is an order-preserving bijection 
  $\theta:\L\rightarrow\M$ with the property that
  $\rank{\theta(M)-\theta(L)} = \rank{M-L}$ for all
  $L<M$ in $\L$. Then there is an invertible operator
  $S$ which maps the range of each $L\in\L$ onto the
  range of $\theta(L)$. Furthermore, $S$ can be taken to be
  an arbitrarily small compact perturbation of a unitary
  operator.
\end{thm}

In addition, the results of \cite{Orr:TrAlIdNeAl}, although
not used directly in the present work, are crucial in the characterization
of stable ideals from \cite{Orr:StIdNeAl}.

\section{preliminaries}

Throughout this paper, $\N$ will denote a continuous nest of projections on
a separable Hilbert space $\H$. See \cite{Davidson:NeAl} for a comprehensive
introduction to nest algebras.

At the heart of our analysis of ideals of continuous nest algebras is a collection
of families of submultiplicative seminorms parameterized by $\N$. These seminorms
have their origins in the earliest work in the field; Ringrose introduced
$i^\pm_N$ in \cite{Ringrose:OnSoAlOp} and showed that the Jacobson radical
is the intersection of the kernels of these seminorms. The following lists
the full set of seminorms we shall need:
 
\begin{defn}
  For $X\in\alg\N$ and $N\in\N$ define
  \begin{align}
    e^+_N(X) & = 
      \begin{cases}
	\lim_{M\downarrow N} \sup_{N<L<M} \enorm{(M-L)X(M-L)} & \text{if $N<I$}, \\
	0                                      & \text{if $N=I$}
      \end{cases} \\
    e^-_N(X) & = 
      \begin{cases}
	\lim_{M\uparrow N} \sup_{N>L>M} \enorm{(L-M)X(L-M)} & \text{if $N>0$}, \\
	0                                      & \text{if $N=0$}
      \end{cases} \\
    i^+_N(X) & = 
      \begin{cases}
	\lim_{M\downarrow N} \norm{(M-N)X(M-N)} & \text{if $N<I$}, \\
	0                                      & \text{if $N=I$}
      \end{cases} \\
    i^-_N(X) & = 
      \begin{cases}
	\lim_{M\uparrow N} \norm{(N-M)X(N-M)} & \text{if $N>0$}, \\
	0                                      & \text{if $N=0$}
      \end{cases} \\
    j_N(X) & = 
      \begin{cases}
	\lim_{M\downarrow N, L\uparrow N} \norm{(M-L)X(M-L)} & \text{if $0<N<I$}, \\
	i_0^+(X)                                      & \text{if $N=0$} \\
	i_I^-(X)                                      & \text{if $N=I$}
      \end{cases}
  \end{align}
  For each $N\in\N$ the seminorms $0$, $e^\pm_N$, $i^\pm_N$, and $j_N$ are called
  the elementary seminorm functions. If the map $(X,N)\mapsto\normn{X}$ is 
  defined so that for each $N\in\N$, $\emptydsf$ takes one
  of these seminorms as a value, then $\emptydsf$ is called a diagonal seminorm
  function.
\end{defn}

\begin{rem}
  For each fixed $N\in\N$, the six elementary seminorms at $N$ form a 
  lattice under pointwise ordering (see Figure~2 of \cite{Orr:StIdNeAl}).
  Thus, with pointwise ordering as $N$ varies, the family of
  all diagonal seminorm functions is a complete lattice. A collection
  $\F$ of diagonal seminorm functions is called a stable family of
  diagonal seminorm functions if it is closed under meets and under
  composition with order automorphisms of $\N$. In other words if
  $\emptydsf^{(i)}\in\F$ ($i=1,2$) and $\theta:\N\rightarrow\N$ is 
  an order automorphism, then
  \[
    \emptydsf^{(1)}\wedge\emptydsf^{(2)}
    \quad\text{and}\quad
    \dsf{\theta(N)}{\;\cdot\;}
  \]
  also belong to $\F$.
\end{rem}

\begin{defn}
  We say that an operator $K\in\alg\N$ is of compact character if
  $(N-M)K(N-M)$ is compact for all $0<M<N<I$ in $\N$. Say that
  an ideal of $\alg\N$ is of compact character if all its elements are
  of compact character.
\end{defn}

In \cite{Orr:StIdNeAl} we showed that $\alg\N$ has exactly eleven stable
ideals of compact character. These ideals (excluding $0$) are listed in 
Figure~1 of \cite{Orr:StIdNeAl}. The following theorem characterizes the 
stable ideals which are not of compact character, and will be the basis for
all our results on stable ideals in the present paper.

\begin{thm}\label{thm:orig-stable-ideal-theorem}
  Let $\I$ be a stable ideal in a continuous nest algebra, $\alg\N$.
  If $\I$ is not of compact character then
  there is a stable family $\F$ of diagonal seminorm functions such that
  $X\in\I$ if and only if, for any $\e>0$, there is a diagonal seminorm
  function $\emptydsf$ in $\F$ such that $\normn{X}<\e$ for all $N\in\N$.
\end{thm}

\section{Characterization of Stable Ideals}\label{sec:main-result}

\begin{defn}
  Let $P_1$ and $P_2$ be two families of intervals of $\N$.
  Say that $P_1$ refines $P_2$, and write $P_1\ge P_2$, if whenever
  $E\in P_1$ there is an interval $F\in P_2$ such that $E\le F$. 
\end{defn}

\begin{defn}
  Let $\O$ be a set of families of intervals of $\N$. 
  We call $\O$ a net of intervals if it is a directed set
  under the ordering of refinement. Call $\O$ a stable net if whenever $\theta$ 
  is an order isomorphism of $\N$ onto itself and $P\in\O$ then the set
  \[
    \theta(P) := \{\theta(E) \st E\in P\}
  \]
  also belongs to $\O$. 
\end{defn}

\begin{prp}\label{prp:nets-specify-ideals}
  If $\O$ is a net of intervals on $\N$ then the set
  \[
    \I := \{X\in\alg\N \st \lim_{P\in\O}\sup_{E\in P} \enorm{EXE} = 0 \}
  \]
  is an ideal of $\alg\N$. If $\O$ is a stable net of intervals, then
  $\I$ is a stable ideal.
\end{prp}

\begin{proof}
  For fixed $X$, the map $P\mapsto\sup_{E\in P} \enorm{EXE}$ is decreasing in $P$,
  and so the limit exists. Since, for any fixed interval $E$, $\enorm{EXE}$ is a submultiplicative seminorm
  on $\alg\N$, then so is $\lim_{P\in\O}\sup_{E\in P} \enorm{EXE}$, and
  $\I$ is its kernel, which is an ideal. 

  If $\ad{S}$ is an automorphism of $\alg\N$ then $S$ induces an order
  isomorphism $\theta:\N\rightarrow\N$ and it is routine to prove that
  \[
    k^{-1} \sup_{E\in P}\enorm{EXE}
      \le \sup_{E\in\theta(P)}\enorm{E(SXS^{-1})E} 
      \le k \sup_{E\in P}\enorm{EXE}
  \]
  where $k=\norm{S}\norm{S^-1}$.
  It is clear from this, and the stability property of $\O$, that $\I$ 
  must be stable under conjugation by $S$.
  
\end{proof}

\begin{rem}
  We shall say that $\I$ is the ideal associated with the net $\O$,
  or that $\I$ arises as the kernel of $\O$.
\end{rem}

\begin{rem}
  We should make a brief remark concerning the essential norm,
  $\enorm{X}$, which is used ubiquitously throughout this paper.
  As a result of \cite[Theorem~5.1]{DavidsonPower:BeApCAl},
  we know that for $T\in\alg\N$,
  \[
    \dist{T}{\K(\H)} = \dist{T}{\K(\H)\cap\alg\N}
  \]
  were $\K(\H)\subseteq\bh$ is the set of all compact operators.
  Thus we shall not make any distinction between the quotient norms
  on $\alg\N/(\K(\H)\cap\alg\N)$ and on $\bh/\K(\H)$, but shall use
  $\enorm{X}$ to denote either, as long as $T\in\alg\N$.
\end{rem}

\begin{eg}
  Let $\O$ consist of the single family, $\{0\}$. Then the associated ideal is
  all of $\alg\N$. Conversely if $\O$ consists of the singleton $\{I\}$ then
  the associated ideal is $\K$, the compact operators in $\alg\N$.
\end{eg}

\begin{eg}
  Let $\O$ consist of the single family, $\{N\in\N \st N < I\}$, of all intervals 
  with lower endpoint equal to $0$. This is a stable net on $\N$ and the 
  associated ideal is the ideal of compact character, $\kminus$. 
  Similarly if $\O$ consists of the singleton $\{M^\perp \st M > 0\}$ then 
  the associated ideal is $\kplus$. Finally, if $\O$ contains the single family
  $\{N-M \st 0 < M < N < I\}$ then the associated ideal is the ideal of all 
  operators of compact character, or $\kplus + \kminus$ 
  (see \cite[Lemma 2.14]{Orr:StIdNeAl}).
\end{eg}

\begin{eg}
  Let $\O$ consist of the set of all singletons $\{N\}$ for $N>0$. This
  is a stable net and the associated ideal is equal to
  \[
    \eo := \{X\in\alg\N\;:\; \inf_{N>0}\norm{NEN}=0\}
  \]
  The transition from the essential norm to the operator norm in this
  example is routine: For if $\enorm{NXN} < \e$ then there is a compact
  $K\in\K$ such that $\norm{NXN + K} < \e$ and we can find a projection
  $0<M<N$ in $\N$ such that $\norm{MKM}<\e$ and so $\norm{MXM} < 2\e$.
  In the same way, the set 
  \[
    \ei := \{X\in\alg\N\;:\; \inf_{N<I}\norm{N^\perp E N^\perp}=0\}
  \]
  also arises as the kernel of a stable net.
\end{eg}

\begin{eg}
  Let $\O$ consist of all finite partitions of $\N$. 
  That is to say, each $P$ in $\O$ is a finite set of
  pairwise orthogonal intervals which sum to the identity. Clearly this is a
  stable net on $\N$, and the associated ideal is the Jacobson radical
  of $\alg\N$ \cite{Ringrose:OnSoAlOp}. The transition from the essential norm to the operator norm
  is similar to the last example. If $\enorm{\sum_{i=1}^m E_iXE_i} < \e$ 
  then there is a compact $K\in\K$ such that 
  $\norm{\sum_{i=1}^m E_iXE_i + K} < \e$. But since $K$ belongs to the 
  Jacobson radical, there is a finite partition $\{F_k\}_{k=1}^n$ refining
  $\{E_i\}$ and such that $\norm{\sum_{i=1}^n F_iKF_i} < \e$. Thus
  $\norm{\sum_{i=1}^n F_iXF_i} < 2\e$.
\end{eg}

\begin{eg}
  Let $\O$ consist of all partitions (finite or infinite) of $\N$.
  Then by \cite{Larson:NeAlSiTr} $\O$ is a net but not a stable net, and
  the associated ideal is Larson's ideal, $\R^\infty_\N$.
\end{eg}

\begin{eg}
  In \cite{Orr:MaIdNeAl} a pseudopartition was defined as a maximal (with respect to
  set inclusion) family of pairwise orthogonal intervals of $\N$. If $\O$ is the
  set of all pseudopartitions of $\N$ then $\O$ is a stable net and 
  the associated ideal is $\I^\infty_N$, which has been shown to be the
  strong radical of $\alg\N$, or the intersection of all the maximal two-side ideals of
  $\alg\N$.
\end{eg}

\begin{prp}\label{prp:intersect}
  Let $\O$ and $\Q$ be two stable nets of intervals on $\N$, associated with
  the stable ideals $\I$ and $\J$ respectively. Then $\I\cap\J$ is the kernel
  of the stable net
  \[
    \O+\Q := \{P\cup Q \st P\in \O, Q\in\Q \}
  \]
\end{prp}

The proof is straightforward, and is left to the reader.

\begin{cor}\label{cor:stable_ideals_of_compact_character}
  All the stable ideals of compact character in $\alg\N$ arise as kernels of 
  stable nets.
\end{cor}

\begin{proof}
  By Theorem~2.16 of \cite{Orr:StIdNeAl}, the stable ideals of compact 
  character are the set of ten ideals listed in Figure~1 of \cite{Orr:StIdNeAl}. 
  Each of these ideals can be expressed as an intersection of the ideals 
  $\kplus+\kminus$, $\kplus$, $\kminus$, $\eo$, $\ei$, and $\K$. 
  The preceding examples have shown that these all arise as the 
  kernels of stable nets, and so the result follows from 
  Proposition~\ref{prp:intersect}.
\end{proof}

\begin{rem}\label{rem:define-induced-stable-net}
  In what follows, let $\I$ be a fixed stable ideal. For each $X\in\I$ and $\e>0$
  define $P_{X,\e}$ to be the set of all intervals $E$ of $\N$ for which
  $\enorm{EXE}<\e$. Finally, let $\O=P_{X,\e} \st X\in\I, \e > 0\}$. 
\end{rem}

\begin{lem}\label{lemma:stable_net}
  With the definitions above, $\O$ is a stable net of intervals on $\N$.
\end{lem}

\begin{proof}
  First we shall show that $\O$ is a directed set. As is thoroughly described in
  \cite{LarsonPitts:IdNeAl,Orr:TrAlIdNeAl}, we can find a projection $F$ in $\N''$
  such that $F$ and $F^\perp$ are both algebraically equivalent to $I$. In other
  words, there are operators $A, B, C, D$ in $\alg\N$ such
  that $AB=F$, $CD=F^\perp$ and $BA=DC=I$, and these operators can be taken to
  be compact perturbations of partial isometries. Now let $P_{X, a}$ and $P_{Y,b}$
  be collections of intervals belonging to $\O$. We shall show that $P_{Z, 1}$ is
  a refinement of $P_{X, a}$ and $P_{Y,b}$, where
  \[
    Z := \frac{1}{a}AXB + \frac{1}{b}CXD
  \]
  To see this, let $E\in P_{Z,1}$ and observe that
  \[
    \begin{split}
      \enorm{EZE} & = \enorm{\frac{1}{a}FE(AXB)EF + \frac{1}{b}F^\perp E(CYD)EF^\perp} \\
                  & = \max \left\{\frac{1}{a}\enorm{EAXBE}, \frac{1}{b}\enorm{ECYDE}\right\}
    \end{split}
  \]
  and this quantity is less than $1$, so that
  \[
    \enorm{EAXBE}<a  \qquad\text{and}\qquad  \enorm{ECYDE}<b 
  \]
  Now since $X=BAXBA$, it follows that $EXE=EBE(AXB)EAE$ and so
  \[
    \enorm{EXE} \le \enorm{EB}\enorm{E(AXB)E}\enorm{AE} \le \enorm{E(AXB)E} < a
  \]
  (with the last inequality following since $A$ and $B$ were taken to be compact 
  perturbations of partial isometries). Thus $E\in P_{X,a}$ and, by the same token,
  $E\in P_{Y, b}$.

  Next, we shall show that $\O$ is stable. In other words, given $P_{X,\e}\in\O$
  and an order isomorphism $\theta$ on $\N$ we must show that $\theta(P_{X,\e})$ is
  also in $\O$. Now by Theorem~\ref{thm:similarity}, find
  an invertible $S$ implementing an automorphism of $\alg\N$ such that $SN=\theta(N)SN$
  for all $N\in\N$. Further take $S$ to be a compact perturbation of a unitary, $S=U+K$. Then
  \[
    \begin{split}
      \enorm{\theta(E)SXS^{-1}\theta(E)} & = \enorm{\theta(E)S(EXE)S^{-1}\theta(E)} \\
					 & \le \enorm{(U+K)EXE(U+K)^{-1}}           \\
					 & = \enorm{EXE}
    \end{split}
  \]
  The same argument, using $S^{-1}$ and $\theta^{-1}$ in place of $S$ and $\theta$,
  yields the reverse inequality and so 
  $\enorm{\theta(E)SXS^{-1}\theta(E)} = \enorm{EXE}$. Thus 
  $\theta(P_{X,\e}) = P_{SXS^{-1}, \e}$, which belongs to $\O$.
\end{proof}

The next definition and the lemma that follows establish a connection between
the collections $P_{T,\e}$ of intervals and diagonal seminorm functions
from Theorem~\ref{thm:orig-stable-ideal-theorem}.

\begin{defn}\label{def:compatible-collections}
  Say that a collection $P$ of intervals of $\N$ is compatible with the diagonal
  seminorm function $\emptydsf$ if, for each $N\in\N$:
  \begin{enumerate}
    \item 
      Whenever $\emptydsf = j_N$ there are projections $G>N>L$ in $\N$
      such that $(G-L)\in P$.
    \item 
      Whenever $\emptydsf=a^+_N\vee a^-_N$ where $a^\pm_N=0, e^\pm_N, i^\pm$ then:
      \begin{enumerate}
	\item
	  If $a^+_N=i^+_N$ then there is a projection $G>N$ in $\N$ 
	  such that $(G-N)\in P$.
	\item
	  If $a^+_N=e^+_N$ then there is a projection $G>N$ in $\N$ 
	  such that $(G-L)\in P$ for all projections $G>L>N$ in $\N$.
	\item
	  If $a^-_N=i^-_N, e^-_N$ then the analogous lower conditions hold.
      \end{enumerate}
  \end{enumerate}
\end{defn}

\begin{lem}\label{lem:compatible-dsf-and-sets-of-intervals}
  Let $T\in\alg\N$ and $a>0$. 
  If $P_{T, \e}$ is compatible with $\emptydsf$ 
  then $\normn{T}\le\e$ for all $N\in\N$.
  Conversely, if $\normn{T}<\e$ for all $N\in\N$ then
  $P_{T, \e}$ is compatible with $\emptydsf$.
\end{lem}

\begin{proof}
  Suppose $P_{T, a}$ is compatible with $\emptydsf$ and
  fix $N\in\N$. Then $A\mapsto\normn{A}$ is one of the 10 possible seminorms listed
  in Figure~2 of \cite{Orr:StIdNeAl}. That is to say, it is of the form
  $a_N(A)\vee b_N(A)$ where $a_N$ is one of $0, e_N^+, i_N^+, j_N$ and
  $b_N$ is one of $0, e_N^-, i_N^-, j_N$. 

  If $a_N=0$ there is nothing to prove. 
  Next, suppose $a_N=e_N^+$.
  Since $P_{T,\e}$ is compatible with $\emptydsf$,
  there is a $G>N$ in $\N$ such that $(G-L)\in P_{T,\e}$
  for all $G>L>N$. Thus
  \[
    a_N(T) \le \sup \{ \enorm{(G-L)T(G-L)} \st L\in\N, N<L<G \} \le \e
  \]
  Next, suppose $a_N=i_N^+$. Since $P_{T,\e}$ is compatible with $\emptydsf$,
  there is a $G>N$ such that $(G-N)\in P_{T,\e}$.
  Thus 
  \[
    \enorm{(G-N)T(G-N)}<\e
  \]
  Pick a compact
  operator $K$ such that $\norm{(G-N)T(G-N) - K}<\e$, and then
  \[
    a_N(T)=i_N^+(T)=i_N^+(T-K)\le\norm{(G-N)T(G-N) - K}<\e
  \]
  The argument for the case $a_N=j_N$ is almost identical.

  Conversely if $\normn{T}<\e$ for all $N\in\N$ then it follows directly
  from the definitions of the seminorm functions that
  we can find intervals of $\N$ with the appropriate
  properties to belong to $P_{T, \e}$.
\end{proof}

\begin{thm}\label{thm:stable_ideals}
  Let $\I$ be a non-zero stable ideal of $\alg\N$. Then there is a stable
  net $\O$ of intervals of $\N$ such that
  \[
    \I = \{ X\in\alg\N \st \lim_{P\in\O}\sup_{E\in\P} \enorm{EXE} = 0 \}
  \]
\end{thm}

\begin{proof}
  By Corollary~\ref{cor:stable_ideals_of_compact_character} we need only consider the case when
  $\I$ is not
  of compact character.
  As above, let $\O$ be the set of collections $P_{T,\e}$ as $T$ ranges over $\I$ and
  $\e$ ranges over all positive values. Since Lemma~\ref{lemma:stable_net} shows $\O$
  is a stable net, it is clear that the limit $\lim_{P\in\O}\sup_{E\in P}\enorm{ETE}$
  exists and is zero for all $T\in\I$.  The main body of the theorem is to establish
  the converse.

  Since $\I$ is a stable ideal, by Theorem~\ref{thm:orig-stable-ideal-theorem} there is a set
  $\F$ of diagonal seminorm functions which specifies $\I$. 
  Now suppose that $\lim_{P\in\O}\sup_{E\in P}\enorm{EXE}=0$ for
  some  $X\in\alg\N$. Given
  $\e>0$, find a $P\in\O$ such that $\enorm{EXE}<\e$ for all $E\in P$. By definition,
  $P=P_{T,a}$ for some $T\in\I$ and $a>0$. By rescaling $T$ we may as well assume
  $P=P_{T,\e}$. 
  Since $T\in\I$, find a diagonal seminorm function $\emptydsf$ in $F$ such that 
  $\normn{T}<\e$ for all $N\in\N$. Thus $\emptydsf$ is compatible with 
  $P_{T,\e}$ and so is compatible with $P_{X, \e}$, which contains
  $P_{T,\e}$. It follows by Lemma~\ref{lem:compatible-dsf-and-sets-of-intervals}
  that $\normn{X}\le\e$ for all $N\in\N$, and so $X\in\I$.
\end{proof}

\section{Cofinal Nets}\label{sec:cofinal-nets}

The net constructed in Theorem ~\ref{thm:stable_ideals} for a general stable ideal
is much larger than any of the natural nets given in the examples of the last section.
This naturally raises the question of when two nets give rise to the same ideal.
Knowing the answer to this question which will also be essential to proving the
quotient norm formula in Theorem~\ref{thm:distance-formula}.
Thus this section will be devoted to establishing the following result:

\begin{thm}\label{thm:cofinal-nets}
  Suppose that $\I_1$ (resp.\ $\I_2$) is the ideal associated with the
  net of intervals $\O_1$ (resp.\ $\O_2$). Then $\I_1\supseteq\I_2$
  if and only if $\O_1$ is cofinal in $\O_2$.
\end{thm}

If $\O_1$ is cofinal in $\O_2$ then clearly for any $T\in\bh$
\[
  \lim_{P\in\O_2}\sup_{E\in P} \enorm{ETE} 
  \ge \lim_{P\in\O_1}\sup_{E\in P} \enorm{ETE} 
\]
and so $\I_1\supseteq\I_2$. To prove the converse, suppose
that $\O_1$ is not cofinal in $\O_2$. This means that there
must be a $Q\in\O_2$ which is not refined by any $P\in\O_1$. The 
strategy in this section will be to suppose that nevertheless
$\I_1\supseteq\I_2$ and derive a contraction. 

Much of the main
part of the argument in this section will rely only on combinatoric
arguments concerned with the ordering of $\reals$. It will be much more convenient to work
directly with real numbers and intervals of real numbers, than with projections
and intervals in $\N$. 
Thus let $\N$ be parameterized as $N(t)$ ($t\in[0,1]$) and for convenience take
$N(t)=I$ for $t\ge 1$ and $N(t)=0$ for $t\le 0$.
Write $E(x,y):= N(y)-N(x)$ and, 
if $P$ is a collection
of open intervals in $\reals$, write $E(P)$ for the
set of intervals $\{E(a,b) \st (a,b)\in P\}$. Write
$\Lambda$ for the set of families $P$ of open intervals
in $\reals$ for which $E(P)$ is in $\O_1$. Observe that
$\Lambda$ is a directed set under refinement, and that
if $\theta:\reals\rightarrow\reals$ is any order isomorphism
then $\theta(P)\in\Lambda$ for all $P\in\Lambda$.

Fix $Q$ as
a collection of open intervals of $\reals$ with the property that
$E(Q)\in\O_2$ and no collection in $\Lambda$ refines it. The proof of 
Theorem~\ref{thm:cofinal-nets} will be established when we see,
in Proposition~\ref{prp:p-refines-q}, that all sufficiently
large $P\in\Lambda$ refine $Q$.

\begin{defn}
  Let $P$ be a collection of open intervals. Say that
  an interval is dominated by $P$ if it is a subset
  of an interval in $P$.
\end{defn}

The next four lemmas establish some basic relationships between
$\Lambda$ and $Q$. Once these facts are in place, no other
reference to operator theory will be used, and the remainder of the section
will be purely combinatoric.
    
\begin{lem}\label{lem:transfer-property-a}
  Given any sequence of pairwise disjoint intervals, none of which
  is dominated by $Q$, we can find a $P\in\Lambda$ which
  also does not dominate any of these intervals.
\end{lem}

\begin{proof}
  Let such a sequence of pairwise disjoint intervals,
  $(x_i, y_i)$, be given.
  For each $i$ choose sequences $x_i<x_j^{(i)}<y_j^{(i)}<y_i$
  with $x_j^{(i)}$ decreasing to $x_i$ and 
  $y_j^{(i)}$ increasing to $y_i$. Let $X_j^{(i)}$ be a finite rank
  partial isometry mapping $N(y_{j+1}^{(i)})-N(y_j^{(i)})$
  into $N(x_{j}^{(i)})-N(x_{j+1}^{(i)})$ and let $X := \sum_{i,j}X_j^{(i)}$.
  Because all the ranges and domains are pairwise orthogonal, the
  sum for $X$ converges weakly to a partial isometry in $\alg\N$.
  Let $(a,b)\in Q$. For each $i$, since no interval of $Q$ contains
  $(x_i, y_i)$, $E(a,b)X_j^{(i)}E(a,b)$ is non-zero for 
  only finitely many $j$.
  Further, since the $(x_i, y_i)$ are pairwise disjoint, no interval
  of $Q$ can meet more than two of them, otherwise it would
  have to contain one of them. Thus in fact $E(a,b)X_j^{(i)}E(a,b)$
  is non-zero for only finitely many values of $i$ and $j$. Hence 
  $E(a,b)XE(a,b)$ is finite rank. Since this holds for
  any $(a,b)\in Q$, we conclude $X\in\I_2$ and so $X\in\I_1$. Thus
  there must be a $P\in\Lambda$ such that 
  \[
    \sup\{\enorm{E(a,b)XE(a,b)} \st (a,b)\in P\} < 1
  \]
  Clearly this is only possible if no interval in $P$ dominates
  any $(x_i, y_i)$.
\end{proof}

\begin{lem}\label{lem:transfer-property-b}
  Given any sequence of pairwise disjoint intervals $(x_i, y_i)$
  with the property that no interval $(x_i, y)$ is dominated
  by $Q$ for any $y>y_i$, then we can find a $P\in\Lambda$ 
  which also does not dominate any interval $(x_i, y)$ with $y>y_i$.
\end{lem}

\begin{proof}
  The proof is very similar to Lemma~\ref{lem:transfer-property-a}.
  Let the sequence of pairwise disjoint intervals 
  $(x_i, y_i)$ be given.
  For each $i$ choose sequences $x_i<x_j^{(i)}<y_i$ and $y_i<y_j^{(i)}$
  with $x_j^{(i)}$ decreasing to $x_i$ and 
  $y_j^{(i)}$ decreasing to $y_i$, and let
  $X_j^{(i)}$ be a finite rank
  partial isometry mapping $N(y_j^{(i)})-N(y_{j+1}^{(i)})$
  into $N(x_{j}^{(i)})-N(x_{j+1}^{(i)})$. Note that this time,
  the intervals $(y_{j+1}^{(i)}, y_j^{(i)})$ need not be pairwise
  disjoint as $i,j$ vary, but we shall simply stipulate that
  the $X_j^{(i)}$ should be chosen with pairwise orthogonal
  initial spaces so that again $X:=\sum_{i,j}X_{i,j}$ converges weakly.
  If $(a,b)\in Q$ then $E(a,b)X_j^{(i)}E(a,b)g$
  is non-zero only if $a<y_i<b$, and this condition can
  be met for at most a single value of $i$. Since 
  $(x_i, y_i)$ is not dominated by $(a,b)$, the compression 
  must be finite rank. Thus 
  $E(a,b)XE(a,b)$ is finite rank, and so 
  $X\in\I_2\subseteq\I_1$. Thus
  there must be a $P\in\Lambda$ such that 
  \[
    \sup\{\enorm{E(a,b)XE(a,b)} \st (a,b)\in P\} < 1,
  \]
  and this is only possible if no interval in $P$ dominates
  any $(x_i, y)$ with $y>y_i$.
\end{proof}
    
The next lemma is almost identical to the last, and the proof is left to the reader.

\begin{lem}\label{lem:transfer-property-c}
  Given any sequence of pairwise disjoint intervals $(x_i, y_i)$
  with the property that no interval $(x_i, x)$ is dominated
  by $Q$ for any $x>x_i$, then we can find a $P\in\Lambda$ 
  which also does not dominate any interval $(x_i, x)$ with $x>x_i$.
\end{lem}

\begin{rem}
  There are obvious analogues to Lemmas~\ref{lem:transfer-property-b}
  and~\ref{lem:transfer-property-c} which deal with the corresponding behavior 
  at the upper endpoints of the intervals $(x_i, y_i)$.
\end{rem}

\begin{lem}\label{lem:p-contained-in-q}
  For all sufficiently large $P\in\Lambda$, $\cup P\subseteq \cup Q$.
\end{lem}

\begin{proof}
  Let $K:=\left(\cup Q\right)^c$. By \cite[Lemma~2.4]{Orr:MaIdNeAl}
  there is an $X\in\alg\N$ such that $\enorm{E(a,b)XE(a,b)}\ge 1$
  if $(a,b)$ intersects $K$ and is zero otherwise (the lemma cited
  claims a sightly weaker result in its statement, but the construction
  used in fact establishes this fact). But thus 
  $E(a,b)XE(a,b)=0$ for $(a,b)\in Q$ and so $X\in\I_2\subseteq\I_1$
  and so there must be a $P\in\Lambda$ with 
  $\sup\{\enorm{E(a,b)XE(a,b)} \st (a,b)\in P\} < 1$. This shows that
  every $(a,b)$ in $P$ must be a subset of $\cup Q$.
  Clearly if $\cup P\subseteq \cup Q$ then the same is true for any
  $P'\in\Lambda$ that refines $P$. 
\end{proof}

\begin{rem}
  Having established some basic relations between 
  the collection $Q$ and at least all sufficiently large
  members of the net $\Lambda$, we shall now develop
  a framework of properties of families of open intervals.
  The next few results will make no assumptions about
  $\Lambda$ or $Q$ and will not use any operator theory.
  We shall return to the context of operator algebras
  with Lemma~\ref{lem:control-endpoints}.
\end{rem}

\begin{defn}
  Given a collection $P$ of open intervals and $x\in\cup P$, write
  \[
    \begin{split}
      L_P(x) := \inf\{a \st x\in(a,b) \text{ for some } (a,b)\in P\} \\
      R_P(x) := \sup\{b \st x\in(a,b) \text{ for some } (a,b)\in P\}
    \end{split}
  \]
  When the context is clear, we shall omit the subscript $P$.
\end{defn}

\begin{rem}\label{rem:props-of-lrfuncts}
  Note that $L(x)$ and $R(x)$ are increasing functions and that 
  for all $x\in\cup P$
  \[
    L(x) < x < R(x).
  \]
\end{rem}

\begin{defn}
  A linked list of intervals is a sequence of intervals $(a_i, b_i)$
  which is indexed by a finite set of integers, or by one of 
  $\integers$, $\integers^+$, or $\integers^-$,
  and has the property
  \[
    a_i < b_{i-1}<a_{i+1}<b_i
  \]
  for all $i$. 
\end{defn}

\begin{rem}
  The union of the the intervals of a linked list is always
  an open interval.
\end{rem}

\begin{defn}
  Let $P$ be a collection of open intervals. 
  Say that an interval is approximately
  dominated by $P$ if it is the union of an increasing
  sequence of intervals each of which is dominated by $P$.
\end{defn}

\begin{rem}
  For any $x\in\cup P$, the intervals $(L(x), x)$ and $(x, R(x))$
  are approximately dominated by $P$.
\end{rem}

\begin{defn}
  Let $P$ be a collection of open intervals with $\cup P=(a,b)$.
  A linked list whose union is $(a,b)$ is called 
  an inner cover of $(a, b)$ if every interval of the list
  is approximately dominated by $P$.
\end{defn}

\begin{defn}
  Let $P$ be a collection of open intervals with $\cup P\subseteq(a,b)$.
  A linked list of intervals whose union is $(a, b)$
  is called an outer cover of $(a, b)$ if every interval of $P$ is contained
  in an interval of the list. 
\end{defn}

Trivially, $(a,b)$ is itself an outer cover for $(a,b)$, so outer covers always
exist. The following lemma shows that inner covers also always exist.

\begin{lem}\label{lem:inner-covers-exist}
  If $P$ is a collection of open intervals and $\cup P = (a,b)$
  then $(a,b)$ has an inner cover with respect to $P$.
\end{lem}

\begin{proof}
  Pick $t_0$ in $(a,b)$ and inductively pick $t_k := R(t_{k-1})$ for
  $k>0$ and $t_k := L(t_{k+1})$ for $k<0$. The sequence so obtained
  is strictly increasing, by Remark~\ref{rem:props-of-lrfuncts}.
  Continue this process for
  as long as these $t_k$ are contained in $(a,b)$.
  By compactness, $t_k$
  increases to $b$ and decreases to $a$, for only finitely many intervals
  are needed to cover any $[a+\e, b-\e]$, and so $t_k$ will be below $a+\e$ or
  above $b-\e$ after finitely many steps.

  Now each $t_k\in(a,b)=\cup P$ and so belongs to an interval of $P$. 
  We can pick pairwise disjoint intervals $(a_k, d_k)$ containing
  $t_k$ and each contained in an interval of $P$. Then choose $b_k$ and $c_k$
  to satisfy 
  \[
    a_k < b_k < t_k < c_k < d_k
  \]
  By construction, each $(c_k, b_{k+1})$ is also contained
  in an interval of $P$. If the sequence of $t_k$ has no greatest
  or smallest element then the sequence of $(a_k, d_k)$ and
  $(c_k, b_{k+1})$ is an inner cover. If the sequence of $t_k$'s
  has a final element, $t_n$, take the final interval of the
  inner cover to be $(c_n, b)$. Likewise if the sequence of $t_k$'s
  has a first element, $t_m$, take the first interval of the inner
  cover to be $(a, b_m)$. These last two intervals need not be 
  dominated by $P$, but since $L(t_m)=a$ and $R(t_n)=b$,
  they clearly are approximately dominated
  by $P$. 
\end{proof}

Our goal is to construct outer covers for $P\in\Lambda$ and
inner covers for intervals of $Q$ which are of compatible order
type (in a sense made precise in Lemma~\ref{lem:order-compatible-bijection}).
As partial steps in that direction, the next two lemmas relate the 
existence of least elements of a cover to 
a property that can be transferred between $P\in\Lambda$ and $Q$ 
using Lemmas~\ref{lem:transfer-property-a} to~\ref{lem:transfer-property-c}.

\begin{lem}\label{lem:inner-cover-with-least-element}
  Let $P$ be a collection of open intervals with $\cup P = (a,b)$.
  The the following are equivalent:
  \begin{enumerate}
    \item There is an $x>a$ with $L(x) = a$.
    \item There is an inner cover of $(a,b)$ with a least element.
  \end{enumerate}
\end{lem}

\begin{proof}
  Suppose that $x>a$ with $L(x) = a$. 
  By Lemma~\ref{lem:inner-covers-exist}
  we know that inner covers for $P$ can be found. Suppose that $(a_i, b_i)$
  is an inner cover with no least element. Then eventually,
  as $i$ decreases, $b_i<x$. 
  Since $L(x)=a$ we know that $(a, x)$ is approximately dominated by $P$
  and so the same is true for $(a, b_i)$.
  It follows that if we delete all intervals to the left of 
  $(a_i, b_i)$ and replace $(a_i, b_i)$ with $(a, b_i)$ then
  we obtain a new inner cover, having a least element.

  On the other hand suppose that there is an inner cover
  with a least element, which we can write $(a, c)$. Since
  $(a,c)$ is the increasing union of intervals which are
  contained in members of $P$
  clearly $L(x) = a$ for all $a<x<c$.
\end{proof}

\begin{lem}\label{lem:outer-cover-with-least-element}
  Let $P$ be a collection of open intervals with $\cup P \subseteq (a,b)$.
  Then the following are equivalent:
  \begin{enumerate}
    \item There is an $x>a$ in $\cup P$ with $L(x) = a$.
    \item Every outer cover of $(a,b)$ has a least element.
  \end{enumerate}
\end{lem}

\begin{proof}
  Suppose that $x>a$ is in $\cup P$ with $L(x)=a$,
  and that $(a_i, b_i)$ is an outer
  cover of $(a, b)$. 
  It follows that $x\in(a_k, b_k)$ for some $k$. 
  If the outer cover has no 
  least element then there is an interval $(a_j, b_j)$ which
  lies to the left of $(a_k, b_k)$ and is not its immediate predecessor.
  Such an interval must satisfy $a<a_j<b_j<x$. However since $L(x)=a$,
  there is an interval $(u,v)$ in $P$ satisfying 
  $a<u<a_j<b_j<x<v$. But no interval of a linked list is a subset of any
  other interval of the list, and so $(u,v)$ cannot be contained in any
  $(a_i, b_i)$, contrary to the property of an outer cover. Thus every
  outer cover must have a least element.

  Suppose there is no $x>a$ in $\cup P$ with $L(x) = a$ and aim
  to construct an outer cover with no least element.
  If it is possible to find a sequence $x_i\in(a,b)\setminus\cup P$
  decreasing to $a$, then we take the cover to be $(x_2, b)$ 
  together with the intervals $(x_{2i+2}, x_{2i-1})$. Otherwise 
  we assume that for all $x$ sufficiently close to $a$, 
  $x\in\cup P$ and $L(x)>a$.
  Then by hypothesis, one can
  inductively choose a sequence $x_i\in\cup P$,
  decreasing to $a$ with the property
  that $x_{i+1}<L(x_i)$ for each $i$. We claim that the sequence
  of intervals $(L(x_{i+1}), x_i)$, together with $(L(x_1), b)$, is an
  outer cover for $P$. Since $L(x_{i+1}) < x_{i+1} < L(x_i) < x_i$,
  the intervals have the correct overlapping property to make them a
  linked list. Now, let an interval $(c,d)$ in $P$ be given and aim to
  show it is contained in an interval of the cover.
  Since $x_i$ decreases to $a$, eventually $x_i<d$, and so let $n$
  be the greatest $i$ with $d\le x_i$. Thus, $x_{n+1}<d\le x_n$ and so 
  $L(x_{n+1})\le c$, so that $(c, d)\subseteq (L(x_{n+1}), x_n)$. 
  If there were no $i$ with $d\le x_i$ then $x_1<d$ and so
  $(c,d)\subseteq(L(x_1), b)$, and the claim is established.
\end{proof}

\begin{rem}
  Clearly there are natural analogues of 
  Lemmas~\ref{lem:inner-cover-with-least-element}
  and~\ref{lem:outer-cover-with-least-element} relating
  the condition $R(x)=b$ to the existence of greatest
  elements in inner and outer covers. 
\end{rem}

\begin{cor}\label{cor:inner-cover-finite}
  Suppose that there are $x,y\in\cap P=(a,b)$ such that
  $L(x)=a$ and $R(y)=b$. Then $(a,b)$ admits a finite
  inner cover.
\end{cor}

\begin{proof}
  By Lemma~\ref{lem:inner-cover-with-least-element} there is an inner
  cover with a least element. If the cover has a greatest
  element we are done, so suppose otherwise. By the dual of 
  Lemma~\ref{lem:inner-cover-with-least-element}, there is another
  inner cover, having a greatest element, $(c, b)$. All
  but finitely many of the intervals from the first cover must
  be contained in $(c, b)$, and so we may form a new, finite
  inner cover consisting of $(c,d)$ together with those
  intervals from the first cover which are not contained in $(c,d)$.
\end{proof}

\begin{cor}\label{cor:outer-cover-ordered-as-z}
  Suppose that $L(x)>a$ and $R(x)<b$ for all $x\in\cup P\subseteq(a,b)$.
  Then there is an outer cover for $(a,b)$ with no greatest or least elements.
\end{cor}

\begin{proof}
  By Lemma~\ref{lem:outer-cover-with-least-element}, there is an outer
  cover $(a_i, b_i)$ with no least element. If this cover
  also has no greatest element we are done, so suppose that
  it has a greatest element, and so without loss let the cover be
  indexed by $i\in\integers^-$. By the dual of
  Lemma~\ref{lem:outer-cover-with-least-element}, there is also
  an outer cover with no greatest element. Likewise
  we shall suppose this cover has a least element and so we
  can list its elements as the intervals $(a'_i, b'_i)$ for $i\in\integers^+$.

  Clearly $a_{-1}<b_{-2}<b_{-1}=b$. Thus $b_{-1}<a'_i<b$
  for all sufficiently large positive values of $i$.
  Select an $i_0$ such that $b_{-1}<a'_{i_0}<b$
  and define $(a_0, b_0)$ to be $(a_{-1}, b'_{i_0})$.
  Finally set $(a_i, b_i) := (a'_{i+i_0}, b'_{i+i_0})$ for
  $i\ge 1$ and we obtain a linked list which is easily seen 
  to be an outer cover.
\end{proof}

\begin{lem}\label{lem:one-element-outer-cover}
  Let $P$ be a collection of open intervals with $\cup P\subseteq (a,b)$.
  The following are equivalent:
  \begin{enumerate}
    \item There is no two-element outer cover of $(a,b)$.
    \item\label{lem:one-element-outer-cover:inner}
      $\cup P = (a,b)$ and there is a one-element inner cover of $(a,b)$.
    \item\label{lem:one-element-outer-cover:interval}
      Every interval $(c,d)$ with $a<c<d<b$ is dominated by $P$.
  \end{enumerate}

\end{lem}

\begin{proof}
  Suppose that $(a,b)$ admits no two-element outer covers and consider
  the sequence of pairs of intervals 
  $(a, b-\frac{1}{n})$ and $(a+\frac{1}{n}, b)$. Since these 
  two intervals can never form
  an outer cover, we must always be able to find an interval 
  $(a_n, b_n)$ in $P$
  which is not contained in either. Such an interval must contain
  $(a+\frac{1}{n}, b-\frac{1}{n})$. 
  This shows that the interval $(a,b)$ is approximately dominated 
  by $P$. Hence the single interval $(a,b)$ is an inner cover for $(a,b)$.

  Conversely, suppose $\cup P = (a,b)$ and $(a,b)$ is an inner cover, and
  suppose for a contradiction that there is an outer cover
  $(a,d), (c, b)$ with $c<d$. Since $(a, b)$ is an inner cover,
  $(a,b)$ is the union of an increasing sequence of intervals dominated by $P$.
  Thus one of these intervals must contain 
  both $c$ and $d$. But this interval is supposed to be a subset of
  a member of $P$, and so be contained in one of
  $(a,d)$ or $(c, b)$, which yields a contradiction.

  The equivalence of 
  items~(\ref{lem:one-element-outer-cover:inner})
  and~(\ref{lem:one-element-outer-cover:interval})
  is an immediate consequence of the definitions.
\end{proof}

\begin{lem}\label{lem:relate-inner-outer-size}
  Let $P$ be a collection of open intervals with $\cup P = (a,b)$
  and suppose there is a sequence $x_1,\ldots, x_n$ in $(a,b)$ satisfying:
  \[
    a=L(x_1),\quad
    R(x_i)=x_{i+1} \; (i=1,\ldots,n-1),\quad
    \text{and}\quad
    R(x_n)=b
  \]
  Then
  \begin{enumerate}
    \item Every outer cover for $(a,b)$ has at most $n+1$ elements.
    \item There exists an outer cover for $(a,b)$ with 
          $\floor{\frac{n}{2}} + 1$ elements.
    \item There exists an inner cover for $(a,b)$ with $n+2$ elements.
  \end{enumerate}
\end{lem}

\begin{proof}
  Suppose if possible that $(a_1, b_1),\ldots, (a_m, b_m)$ is an 
  $m$-element outer cover where $m>n+1$.
  Since $L(x_1)=a$, there is an interval $(c,d)\in P$ containing $x_1$
  with $c<b_2$. The only interval of the outer cover that can contain this
  is $(a_1, b_1)$, and so $x_1\in(a_1, b_1)$. 

  Having now established that $x_1\le b_1$, suppose for induction that
  $x_k\le b_k$. If it were possible that $b_{k+1}<x_{k+1}$ then, since
  $x_{k+1} = R(x_k)$, there would be an interval 
  $(e, f)$ of $P$ containing $x_k$ with $b_{k+1}<f$. But $(e, f)$
  must be contained in an interval $(a_i, b_i)$ satisfying
  \[
    a_i \le e < x_k \le b_k < b_{k+1} < f \le b_i
  \]
  Thus $(a_i, b_i)$ 
  contains points greater than $b_{k+1}$. By the ordering
  property of linked lists, this implies that $i\ge k+2$.
  Similarly we see $(a_i, b_i)$contains points smaller than $b_k$ which,
  implies $i\le k+1$. From this contradiction
  we see by induction that $x_i\le b_i$ for all
  $i=1,2,\ldots, n$. But since $R(x_n)=b$ there is an interval
  $(c, d)$ in $P$ containing $x_n$ with $b_{m-1}<d$. 
  The only interval of the outer cover than can contain $(c,d)$ is $(a_m, b_m)$.
  But since we have also seen $x_n\le b_n$, this means that $(a_m, b_m)$
  must meet an $(a_i, b_i)$ with $i\le n<m-1$ or, in other words, $i\le m-2$. 
  This is impossible for any linked list, and so $m\le n+1$.

  Next we shall construct an outer cover having $m := \floor{\frac{n}{2}}+1$
  elements. Consider the intervals $E_i := (x_{2i-3}, x_{2i})$ for
  $i=1,2,\ldots,m$ (where we define $x_i$ to be $a$ for $i\le 0$ and to be
  $b$ for $i\ge n$). One readily sees that this sequence is 
  a linked list. To see that it is an outer cover, let $(c, d)\in P$ be given
  and let $i_0$ be the smallest $i$ for which $c<x_i$. If $i_0=n+1$ then
  $x_n\le c<d\le x_{n+1}=b$ and so $(c,d)\subseteq E_m$. Otherwise, 
  $1\le i_0 \le n$, and either $d\le x_{i_0}$, or else $c<x_{i_0}<d$, which
  implies $d\le R(x_{i_0})=x_{i_0+1}$. In either case, 
  $(c,d)\subseteq(x_{i_0-1}, x_{i_0+1})$. Every such interval is contained in
  an $E_i$, so the $E_i$ are an outer cover.

  Finally let us construct an inner cover with $n+2$ elements.
  First take the two intervals $(a, x_1)$ and $(x_n, b)$.
  Next, since $x_n\in(a,b)=\cup P$, find an interval 
  $(a_n, b_n)$ in $P$ that contains $x_n$.
  If necessary, adjust $a_n$, $b_n$ so that
  $x_{n-1}<a_n<x_n<b_n<b$ and so that $(a_n, b_n)$ is still
  dominated by $P$.
  Then start at $x_{n-1}$ and work backwards, choosing intervals
  $(a_i, b_i)$ dominated by $P$, that contain $x_i$ and have the property
  $x_{i-1}<a_i<x_i$ and
  $a_{i+1}<b_i\le x_{i+1}$. Stepping backwards, we terminate
  at $i=1$, having found a list of $2+1+(n-1)=n+2$ intervals.
\end{proof}

\begin{cor}~\label{cor:outer-cover-from-pw-disjoint}
  Let $P$ be a collection of open intervals with $\cup P\subseteq (a,b)$.
  If there is a sequence of pairwise disjoint nonempty subintervals $(a_i, b_i)$ 
  $i=1,2,\ldots,k$ of $(a,b)$ none of which are dominated by $P$, then
  $(a,b)$ has an outer cover of size $\ceil{k/2}$.
\end{cor}

\begin{proof}
  First add intervals to $P$ so that in fact $\cup P(a,b)$, without
  in the process changing the fact that no interval $(a_i, b_i)$ is dominated by 
  $P$. This is easy to accomplish by, for example, ensuring that every interval
  added is shorter than all of the $(a_i, b_i)$. Now if there is no $x\in(a,b)$
  for which $L(x)=a$ then by Lemma~\ref{lem:outer-cover-with-least-element} there
  is an infinite outer cover. Likewise there is an infinite outer cover if 
  there is no $x\in(a,b)$ for which $R(x)=b$, and so we may suppose
  that $L(x_1)=a$ for some $x_1>a$ and that $R(x)=b$ for $x$ sufficiently close to $b$.
  Without loss,
  assume that the $(a_i, b_i)$ have been indexed so that the $a_i$ (and $b_i$) are strictly
  increasing. Since we must have $a<b_1$, we may also suppose $x_1<b_1$. 

  Now recursively define $x_{i+1} := R(x_i)$ for as long as the sequence lies in $(a,b)$. 
  This sequence must terminate, for it is strictly increasing and cannot have a limit
  point inside $(a,b)$. Thus eventually $x_i$ must increase to a value at which $R(x_i)=b$,
  at which point the sequence terminates. Suppose that the sequence has $n$ terms.

  Observe that whenever $x_i \le b_i$
  then $x_{i+1}=R(x_i)$ cannot be greater than $b_{i+1}$, otherwise $(a_{i+1}, b_{i+1})$
  would be dominated by $P$. Thus inductively $x_i \le b_i$
  for $1 \le i \le \min\{k, n\}$. But if $k>n+1$ then this shows that
  $x_n\le a_{n+1}<b_{n+1}\le a_{n+2}<b$, and, since $R(x_n)=b$, this shows
  $(a_{n+1},b_{n+1})$. It follows from this contradiction that in fact
  $n\ge k-1$. By Lemma~\ref{lem:relate-inner-outer-size}, there exists an
  open cover of at least $\floor{\frac{k-1}{2}}+1= \ceil{k/2}$ terms.
\end{proof}

\begin{lem}\label{lem:order-compatible-bijection}
  Suppose that $E_i$ ($i\in I$) and $F_j$ ($j\in J$) are two
  linked lists and that $I\supseteq J$. 
  Then there is an increasing bijection
  $\theta:\reals\rightarrow\reals$ which maps each interval
  $E_i$ into an interval $F_j$.
\end{lem}

\begin{proof}
  Both $I$ and $J$ are subsets of $\integers$.
  For each $j\in J$, let $C_j$ be the set of 
  $i\in I$ which are closer to $j$ than to any
  other element of $J$. (In the event $i$
  is equidistant from two elements of 
  $J$, assign it to the smaller of the two.)
  The $C_j$ partition $I$ into ranges of consecutive
  numbers.
  Each of the sets $G_j := \bigcup_{i\in C_j}E_i$
  is an interval, and the collection of intervals so formed
  is a linked list.
  
  Write $G_j = (a_j, b_j)$ and $F_j = (a'_j, b'_j)$.
  By the overlapping property of linked lists
  we can define an order preserving map that
  takes each $a_j\mapsto a'_j$ and each $b_j\mapsto b'_j$. 
  This correspondence can be extended to a piecewise linear
  bijection of $\reals\rightarrow\reals$ mapping each $G_j$
  to $F_j$.
\end{proof}

\begin{defn}
  Let $P$ and $Q$ be two collections of open subintervals of $(a,b)$.
  We shall say that $P$ and $Q$ are order compatible if we can
  find an outer cover $E_i$ ($i\in I$) of $(a,b)$ with respect to
  $P$ and an inner cover $F_j$ ($j\in J$) of $(a,b)$ with respect
  to $Q$ that satisfy $I\supseteq J$. 
\end{defn}

We shall now return to the analysis of the two fixed stable ideals
$\I_1$ and $\I_2$, with the associated net $\Lambda$ and fixed
family of open intervals $Q$. Recall in Lemma~\ref{lem:p-contained-in-q} 
we saw that $\cup P\subseteq\cup Q$ for all sufficiently
large $P\in\Lambda$. We shall fix the following notation for the
remainder of this section. Write $(a_i, b_i)$ ($i\in\naturals$)
for a fixed enumeration
of the connected components of $\cup Q$. If $P$ is any collection
of open intervals satisfying $\cup P\subseteq\cup Q$ we shall write
$\cpt{P}{i}$ for the set of all members of $P$ which are subsets
of $(a_i, b_i)$. 

Our main goal in the remainder of this section
is to show, in Proposition~\ref{prp:p-refines-q}, that all
sufficiently large $P\in\Lambda$ refine $Q$. To do this, we
shall find a $P$ with the property that $\cup P\subseteq\cup Q$
and use the machinery of the last few lemmas to show each $\cpt{P}{i}$
refines $\cpt{Q}{i}$. The main intermediate result is to show
in Lemma~\ref{lem:make-order-compatible} that we can ensure
all the $\cpt{P}{i}$ are order compatible with the corresponding $\cpt{Q}{i}$.
The next couple of lemmas are needed to establish this result.

For any collection of open sets write 
\[
  A(P) := \{i \st L_{\cpt{P}{i}}(x) > a_i \text{ for all } x>a_i \} 
\]
and
\[
  B(P) := \{i \st R_{\cpt{P}{i}}(x) < b_i \text{ for all } x<b_i \} 
\]

\begin{lem}\label{lem:control-endpoints}
  For all sufficiently large $P\in\Lambda$, $A(P)\supseteq A(Q)$ and $B(P)\supseteq B(Q)$.
\end{lem}

\begin{proof}
  We shall prove only the inclusion for $A(P)$. The result for
  the $B(P)$ will follow by dual arguments.
  If $i\in A(Q)$ then $L_{\cpt{Q}{i}}(x)>a_i$ for all
  $x>a_i$. Thus for each $i\in A(Q)$ we can inductively construct sequences 
  $x^{(i)}_k$, $y^{(i)}_k$ 
  in $(a_i, b_i)$ which decrease to $a_i$ and satisfy
  $x^{(i)}_{k+1} < y^{(i)}_{k+1} < L_{\cpt{Q}{i}}(x^{(i)}_k)$ for all $k$. 
  The intervals $I_{i,j} := (y^{(i)}_{k+1}, y^{(i)}_k)$ cannot be dominated by $Q$,
  and are pairwise disjoint as $i$ and $k$ run over all possible values. 
  Thus by Lemma~\ref{lem:transfer-property-a}
  there must be a $P_0$ in $\Lambda$ which does not dominate any $I_{i,j}$.
  Likewise, if $P\ge P_0$, then $P$ does not dominate any $I_{i,j}$.
  But if there were any $i\in A(Q)$ for which we could find an $x>a_i$ 
  with $L_{\cpt{P}{i}}(x)=a_i$, then
  $P$ would dominate all $I_{i,j}$ contained in $(a_i, x)$. 
  It follows that when $i\in A(Q)$ 
  then for any $x>a_i$, $L_{\cpt{P}{i}}(x)>a_i$. Thus $i\in A(P)$.
\end{proof}

\begin{lem}\label{lem:double-outer-covers}
  Suppose we are given some $P_0\in\Lambda$ and a collection $C$ of indices
  such that, for each $i\in C$, $(a_i, b_i)$ has a finite outer cover 
  of size $n_i>1$ with respect to $\cpt{P_0}{i}$.
  Then there is a $P_1$ such that, for all $P\ge P_1$ and all $i\in C$,
  $(a_i, b_i)$ has an outer cover of size $2n_i-1$ with respect to $\cpt{P}{i}$.
\end{lem}

\begin{proof}
  Fix $i$ and let $(x_k, y_k)$ ($k=1,2,\ldots, n_i)$ be an outer
  cover for of $(a_i, b_i)$ with respect to $\cpt{P_0}{i}$. 
  Since these intervals are a linked list, the endpoints
  satisfy
  \begin{equation}\label{lem:double-outer-covers:xyrels}
    x_k < y_{k-1} < x_{k+1} < y_k\text{ for $2\le k\le n_i-1$}
    \text{ and }
    x_1=a_i,\, y_{n_i}=b_i
  \end{equation}
  Thus we can easily pick 
  $x'_k$ ($2\le k\le n_i$) and $y'_k$ ($1\le k\le n_i-1$) to satisfy
  \begin{equation}\label{lem:double-outer-covers:xyprimerels}
    y_{k-1} < x'_k < y'_{k-1} < x_{k+1}\text{ for $2\le k\le n_i-1$} 
  \end{equation}
  and
  \[
    y_{n_i-1} < x'_{n_i} < y'_{n_i-1} < b_i
  \]
  Define $x'_1:=a_i$ and $y'_{n_i}:=b_i$.
  Then take a continuous increasing bijection that maps each $x_k$ 
  to $x'_k$ 
  and each $y_k$ to $y'_k$ for $1\le k\le n_i$. 
  We can construct a single function
  $\reals\rightarrow\reals$ which accomplishes the corresponding mapping
  on each $(a_i, b_i)$. Let $P'_0$ be the image of $P_0$ under this transformation.
  Since $\Lambda$ is a stable net, $P'_0\in\Lambda$. Let $P_1$ be a collection
  in $\Lambda$ that refines both $P_0$ and $P'_0$ and let $P$ be an arbitrary
  collection in $\Lambda$ that refines $P_1$. We shall show that each $\cpt{P}{i}$
  has an outer cover of size $2n_i-1$.

  If $E$ is an interval in $\cpt{P}{i}$ then it is contained in an interval of
  $P_0$ and in an interval of $P'_0$. Since the $(x_j,y_j)$ are an outer cover
  for $(a_i,b_i)$ with respect to $\cpt{P_0}{i}$ 
  and, correspondingly, the $(x'_k, y'_k)$ are an outer cover
  for $(a_i,b_i)$ with respect to $\cpt{{P'_0}}{i}$, it follows that for some $j$ and $k$,
  \[
    E\subseteq(x_j,y_j)\cap(x'_k, y'_k)
  \]
  However by~(\ref{lem:double-outer-covers:xyrels})
  and~(\ref{lem:double-outer-covers:xyprimerels}), the only way this
  intersection can be non-empty is for it to equal one of
  \[
    (x_j, y'_{j-1})\quad\text{ or }\quad(x'_j, y_j)
  \]
  It is routine to check 
  (again by~(\ref{lem:double-outer-covers:xyrels})
  and~(\ref{lem:double-outer-covers:xyprimerels}))
  that the collection of intervals $(x'_l, y_l)$ ($1\le l\le n_i$)
  together with $(x_{l+1}, y'_l)$ ($1\le l \le n_i-1$) is a linked list
  of the correct length, and the result follows.
\end{proof}

\begin{rem}
  Note that $n_i\ge 2$ in the hypotheses of 
  Lemma~\ref{lem:double-outer-covers}, and so
  in the conclusion $2n_i-1\ge\frac{3}{2}n_i$. In the sequel
  we shall need to apply Lemma~\ref{lem:double-outer-covers}
  repeatedly, and this lower bound on the size of the
  outer cover will be easier to iterate.
\end{rem}

\begin{lem}\label{lem:make-order-compatible}
  There is a $P_0$ in $\Lambda$ such that for all $P\ge P_0$ in $\Lambda$ 
  and all $i$, $\cpt{P}{i}$ is order compatible with $\cpt{Q}{i}$. 
\end{lem}

\begin{proof}
  By Lemma~\ref{lem:control-endpoints}, for all sufficiently large $P$,
  $A(Q)\subseteq A(P)$ and $B(Q)\subseteq B(P)$. We shall first show that
  for such a $P$, $\cpt{P}{i}$ is order compatible with $\cpt{Q}{i}$
  for all $i\in A(P)\cup B(P)$. For if in fact $i\in A(P)\cap B(P)$
  then by Corollary~\ref{cor:outer-cover-ordered-as-z} there is
  an outer cover ordered as $\integers$, which clearly is order compatible
  with any inner cover for $\cpt{Q}{i}$. On the other hand if
  $i\in A(P)\setminus B(P)$ then by Lemma~\ref{lem:outer-cover-with-least-element}
  there exists an outer cover with respect to  $\cpt{P}{i}$ which is ordered as $\integers^-$. 
  Correspondingly, since $i\not\in B(Q)$, it follows by
  Lemma~\ref{lem:inner-cover-with-least-element} that there must exist an inner
  cover with respect to
  $\cpt{Q}{i}$ which is either finite or ordered as $\integers^-$. In either
  case $\cpt{P}{i}$ is order compatible with $\cpt{Q}{i}$. Finally, if
  $i\in B(P)\setminus A(P)$, compatibility follows by dual arguments.

  Thus we can restrict attention
  to the case $i\not\in A(P)\cup B(P)$. In such a case, by 
  Corollary~\ref{cor:inner-cover-finite}, since $i\not\in L(Q)\cup R(Q)$,
  $(a_i, b_i)$ has a finite inner cover with respect to $\cpt{Q}{i}$. 
  Suppose in each such case we have picked an inner cover of least cardinality, $m_i$.
  Because there is a finite inner cover, we can find a sequence
  \[
    a_i=L_{\cpt{Q}{i}}(x_1) < x_1
      < R_{\cpt{Q}{i}}(x_1) = x_2
      < R_{\cpt{Q}{i}}(x_2) = x_3 
      < \cdots 
      < R_{\cpt{Q}{i}}(x_n) = b
  \]
  and by Lemma~\ref{lem:relate-inner-outer-size} together with
  the minimality of $m_i$, we have $n\ge m_i-2$.
  Thus we can choose $k_i := \floor{(m_i-2)/2}\ge (m_i-3)/2$ 
  pairwise disjoint subintervals of $(a_i, b_i)$,
  each of which contains both $x_{2j-1}$ and $x_{2j}$ for some $j$.
  None of these intervals can be dominated by
  $\cpt{Q}{i}$. Write $(a^{(i)}_j, b^{(i)}_j)$ 
  ($j=1,2,\ldots, k_i$) for these intervals.
  Thus, by Lemma~\ref{lem:transfer-property-a}, for all sufficiently
  large $P\in\Lambda$, no interval of $\cpt{P}{i}$ contains any 
  $(a^{(i)}_j, b^{(i)}_j)$. This shows, by Corollary~\ref{cor:outer-cover-from-pw-disjoint},
  that there is an outer cover of $(a_i, b_i)$ with respect to $\cpt{P}{i}$
  that has at least $n_i := \ceil{k_i/2}\ge (m_i-3)/4$ elements.

  Let $C$ be the set of $i\not\in A(P)\cup B(P)$ for which
  $m_i>15$. Then for all sufficiently large $P$, 
  $\cpt{P}{i}$ has a finite outer cover of size at least $n_i$ for all $i\in C$, and
  $n_i \ge (m_i-3)/4 > 1$. 
  Thus we can apply Lemma~\ref{lem:double-outer-covers}
  four times, and conclude that for all sufficiently large $P$,
  each $\cpt{P}{i}$ has an outer cover of size at least
  \[
    \left(\frac{3}{2}\right)^4 n_i 
      > 5n_i
      \ge \frac{5m_i-15}{4}
      > m_i
  \]
  It remains to deal with those $i$ for which $m_i\le 15$.

  Consider the set $C$ of all $i$ for which $1<m_i\le 15$.
  By 
  Lemma~\ref{lem:one-element-outer-cover}, 
  we can find $a_i<c_i<d_i<b_i$ such that
  the interval $(c_i, d_i)$ is not dominated by 
  $\cpt{Q}{i}$. By Lemma~\ref{lem:transfer-property-a}, 
  for all sufficiently
  large $P$, no $\cpt{P}{i}$ with $i\in C$ contains $(c_i, d_i)$.
  Thus, again by Lemma~\ref{lem:one-element-outer-cover}, all $(a_i, b_i)$
  with $i\in C$ have two-element outer covers with respect to $\cpt{P}{i}$. 
  Applying Lemma~\ref{lem:double-outer-covers} repeatedly
  five times, we conclude that for all sufficiently large
  $P$, we can find outer covers for all $\cpt{P}{i}$ ($i\in C$) of size at least
  $(\frac{3}{2})^5\times 2 > 15$.

  For all remaining $i$, $m_i=1$ and since every $\cpt{P}{i}$ trivially has
  an outer cover of length 1, we are done. 
\end{proof}

\begin{prp}\label{prp:p-refines-q}
  For all sufficiently large $P\in\Lambda$, $P$ refines $Q$.
\end{prp}

\begin{proof}
  By Lemma~\ref{lem:make-order-compatible} we can find a $P_0\in\Lambda$
  such that each $\cpt{P_0}{i}$ is order compatible with $\cpt{Q_0}{i}$.
  By Lemma~\ref{lem:order-compatible-bijection} we can find order preserving
  bijections defined on each $(a_i, b_i)$ which map the intervals
  of an outer cover of $(a_i, b_i)$ with respect to $\cpt{P_0}{i}$
  into the intervals of an inner cover of $(a_i, b_i)$ with respect 
  to $\cpt{Q}{i}$. We can patch these maps together and
  extend to a single map $\theta:\reals\rightarrow\reals$ which maps
  an outer cover on each $(a_i, b_i)$ into an inner cover. 

  Next, transform $P_0$ to $P_1 := \theta(P_0)$, which belongs to 
  $\Lambda$ because of the stability property of $\Lambda$. Now
  pick an arbitrary interval $E$ which is dominated by $P_1$, 
  and aim to show $E\in Q$. Since by definition the $(a_i, b_i)$
  are the connected components of $\cup Q$, we can find $i$ such that
  $E\subseteq(a_i, b_i)$. Thus $E\subseteq\theta(E')$ for some
  $E'\in \cpt{P_0}{i}$. There is an outer cover of $(a_i, b_i)$
  with respect to $\cpt{P_0}{i}$ which $\theta$ maps into an inner cover
  of $(a_i, b_i)$ with respect to $\cpt{Q}{i}$. Since $E'$ is a subset of
  an interval of the outer cover, $\theta(E')$ is a subset of an interval
  of the inner cover, and so the same is true for $E$. If it so happens
  that $E=(x,y)$ where $a_i<x$ and $y<b_i$ then $E$ must in fact be
  dominated by $Q$, and we are done. Thus we shall consider the
  case when $E=(a_i, x)$. The case of $E=(x, b_i)$ is analogous.

  Let $C$ be the set of $i$ for which no interval $(a_i, x)$ belongs
  to $Q$ for any $x>a_i$. Applying Lemma~\ref{lem:transfer-property-c} 
  to the intervals $(a_i, b_i)$ shows that we can find $P_2\ge P_1$
  such that, for $i\in C$, no interval $(a_i, x)$ is dominated by $P_1$
  for any $x>a_i$. Thus for $E=(a_i, x)$, provided at least that $E$ 
  is dominated by $P_2$, then we know $i\not\in C$.

  For $i\not\in C$, define $c_i := \sup\{x \st (a_i, x)\in Q\}$
  and let $A:=\{i\not\in C \st (a_i, c_i) \not\in Q\}$ 
  and $B:=\{i\notin C \st (a_i, c_i) \in Q\}$.
  Apply Lemma~\ref{lem:transfer-property-a} 
  to the intervals $(a_i, c_i)$ for $i\in A$, and apply
  Lemma~\ref{lem:transfer-property-b} 
  to the intervals $(a_i, c_i)$ for $i\in B$, and so conclude there
  is a $P_3\ge P_2$ which dominates no $(a, c_i)$ (for $i\in A$) or
  $(a_i, x)$ (for $i\in B$ and $c>c_i$). Thus, provided $E=(a_i, x)$ 
  is dominated by $P_3$, we know that $i$ must be in $A$ or $B$,
  and that in either of these cases, $E$ is in $Q$.

  After applying a similar argument to deal with the case $E=(x, b_i)$,
  we finally obtain  $P_4\ge P_3$ with the property that every
  interval dominated by $P_4$ must be dominated by $Q$, and we are done.
\end{proof}

With Proposition~\ref{prp:p-refines-q} we have established the final step of the proof of
Theorem~\ref{thm:cofinal-nets}. However in Proposition~\ref{prp:nets-specify-ideals}
we saw that nets which may not be stable still give rise to
ideals (though not stable ideals), and so it is natural to
ask whether the conclusion of Theorem~\ref{thm:cofinal-nets}
holds without the assumption of stability. The following 
example\footnote{
  Thanks are due to David Pitts for
  suggesting that this example should be included.
} shows that it does not. 

\begin{eg}
  Let $\N\setminus\{0,I\}$ be parameterized by $\reals$ with the
  strongly continuous mapping $t\mapsto N_t$ and, for any 
  $S\subseteq\reals$, define
  \[
    P_S := \{N_{s+1} - N_s \st s\in S\}
  \]
  Let $\O_1$ be the set of all $P_S$ where $S\subseteq\rationals$ and
  $\rationals\setminus S$ is a $G_\delta$ set, and let $\O_2$ be the 
  singleton $\{ P_{\reals\setminus\rationals} \}$. One readily verifies
  that $\O_1$ and $\O_2$ are nets (although not stable nets) of intervals
  and that $\O_1$ is not cofinal in $\O_2$ (or vice versa, for that matter).
  Neverthless we shall show that if $\I_1$ and $\I_2$ are the ideals
  induced by $\O_1$ and $\O_2$ respectively, then $\I_1\supseteq\I_2$.

  To see this, suppose that $X\not\in\I_1$, and aim to show
  $X\not\in\I_2$. Since $X\not\in\I_1$, there is an $\e_0>0$ 
  such that for every $G_\delta$ subset $L$ of $\rationals$,
  \[
    \sup\left\{\enorm{(N_{s+1} - N_s)X(N_{s+1} - N_s)} \st s\in\rationals\setminus L\right\}
      \ge \e_0
  \]
  Now consider the set $T$ of all $x\in\reals$ for which
  \[
    \inf_{t>0}\enorm{(N_{x+t} - N_x)X(N_{x+1} - N_{x+1-t})} 
      \ge \frac{\e_0}{3}
  \]
  If $T$ contains any irrationals then $\enorm{(N_{s+1} - N_s)X(N_{s+1} - N_s)}$
  must be at least $\e_0/3$ for some $s\in\reals\setminus\rationals$,
  and so $X\not\in\I_2$. Thus for the remainder of this argument
  we can assume $T\subseteq\rationals$.

  Recall we are working over a separable Hilbert space, and so let
  $F_j$ be a countable norm dense sequence in the set of compact operators.
  It follows that $x\in T$ if and only if
  \[
    \norm{(N_{x+\frac{1}{i}} - N_x)X(N_{x+1} - N_{x+1-\frac{1}{i}})-F_j} 
      > \frac{\e_0}{3} - \frac{1}{k}
  \]
  for all $i,j,k\in\naturals$. By strong upper continuity of the norm,
  the set of $x$ satisfying this last inequality is an open set, 
  and so $T$ is a $G_\delta$ subset of $\rationals$. But of course,
  by the Baire Category Theorem, $\rationals$ is not itself a $G_\delta$ set
  and so $\rationals\setminus T$ is non-empty and we can find an
  $s\in\rationals\setminus T$ such that 
  $\enorm{(N_{s+1} - N_s)X(N_{s+1} - N_s)}\ge2\e_0/3$. However since
  $s\not\in T$, for all sufficiently small $\eta>0$,
  \[
    \enorm{(N_{s+\eta} - N_s)X(N_{s+1} - N_{s+1-\eta})} < \frac{\e_0}{3}
  \]
  and so we shall fix on a small irrational value of $\eta$ for which this holds.
  From the last two inequalities it follows that 
  \[
    \enorm{
      (N_{s+1} - N_s)X(N_{s+1} - N_s) 
        - (N_{s+\eta} - N_s)X(N_{s+1} - N_{s+1-\eta})
    }
    \ge \frac{\e_0}{3}
  \]
  and therefore at least one of the following must hold:
  \[
    \enorm{(N_{s+1-\eta}-N_{s-\eta})X(N_{s+1-\eta}-N_{s-\eta})}
     \ge \frac{\e_0}{6}
  \]
  or
  \[
    \enorm{(N_{s+1+\eta}-N_{s+\eta})X(N_{s+1+\eta}-N_{s+\eta})}
     \ge \frac{\e_0}{6}
  \]
  However since $s\pm\eta$ is irrational this means
  \[
    \sup\{\enorm{EXE} \st E\in P_{\reals\setminus\rationals}\} \ge \frac{\e_0}{6}
  \]
  and so, since $\O_2$ is a singleton, the limit over $\O_2$ is non-zero,
  and $X\not\in\I_2$.
\end{eg}

\section{Quotient Norms}\label{sec:quotient-norms}

The main result of this section will be Theorem~\ref{thm:distance-formula},
which establishes a formula for the quotient form of $\alg\N$ by
a stable ideal. We shall first, in Corollaries~\ref{cor:kminus-quotient}
and~\ref{cor:kplus-quotient} establish a version of the quotient
norm formula for two ideals of compact character, $\kplus$ and $\kminus$.
The following lemma, which is derived from a theorem of
Axler, Berg, Jewell, and Shields 
\cite[Theorem~2]{AxlerBergJewellShields:ApCoOpSpH} 
is needed for these formulas.

\begin{lem}\label{lem:norm-estimate}
  Let $X\in\bh$ and let $P$ be a projection such that
  \[
    \enorm{X}<a\qquad\text{and}\qquad\norm{XP^\perp}<a.
  \]
  Suppose $E_n=E_nP$ is a sequence of operators converging strong-* to $P$
  and satisfying $\lim_n\norm{I-E_n}= 1$.
  Given any $\e>0$ there is an $n_0$ such that 
  $\norm{X(I-E_n)}\le a+\e$ for all $n\ge n_0$.
\end{lem}

\begin{proof}
  Rescaling as necessary, we shall assume $a=1$
  and suppose for a contradiction there is an $\e_0>0$ such that 
  $\norm{X(I-E_n)}>1+\e_0$ 
  for infinitely many $n$. For each such $n$ we pick a unit vector $t_n$
  satisfying $\norm{X(I-E_n)t_n}>1+\e_0$ and set $x_n := (I-E_n)t_n$.
  Passing to a subsequence we may assume that this inequality holds
  for all $n$ and, further, that $x_n$ is weakly convergent to a limit $x$. Write
  $x_n=x+e_n$ where $\wlim e_n = 0$. We claim that $x= P^\perp x$.
  To see this, observe that for any fixed $m$ and $y\in\H$, as $n\rightarrow \infty$
  \[
    |\innerprod{E_m(I-E_n)t_n}{y}|
      \le \norm{(I-E^*_n)E^*_m y} 
      \longrightarrow \norm{(I-P)E^*_m y} = 0
  \]
  and thus $E_m x = \wlim_{n\rightarrow\infty} E_m x_n = 0$ for all $m$. Taking the limit as 
  $m\rightarrow\infty$, our claim follows.  

  Since $\enorm{X}<1$, find a finite rank projection $F$ such that $\norm{F^\perp X}<1$.
  Without loss we may assume $F$ includes $Xx$ in its range. Finally, pick $n_0$ such
  that
  \[
    \norm{x_n} < 1 + \frac{\e_0}{4},\qquad
    |\innerprod{e_n}{x}| < \frac{\e_0}{8}
    \qquad\text{and}\qquad
    \norm{FXe_n}< \frac{\e_0}{2}
  \]
  for all $n\ge n_0$. It now follows that for any $n\ge n_0$,
  \[
    \begin{split}
      \norm{Xx_n} & =   \norm{FXx + FXe_n + F^\perp Xx + F^\perp Xe_n} \\
                  & \le \norm{FXx + F^\perp Xe_n} + \frac{\e_0}{2}
    \end{split}
  \]
  and
  \[
    \begin{split}
      \norm{FXx + F^\perp Xe_n}^2 
        & = \norm{FXx}^2 + \norm{F^\perp Xe_n}^2             \\
        & = \norm{FXP^\perp x}^2 + \norm{F^\perp Xe_n}^2     \\
        & < \norm{x}^2 + \norm{e_n}^2 
            \qquad\text{(since $\norm{XP^\perp}<1$ and $\norm{F^\perp X}<1$)}    \\
        & \le \norm{x_n}^2 + 2|\innerprod{e_n}{x}|           \\
        & < 1 + \frac{\e_0}{2} < \left(1 + \frac{\e_0}{2}\right)^2
    \end{split}
  \]
  Thus, $\norm{Xx_n}<1 + \e_0$, contradicting our hypothesis.
\end{proof}

\begin{cor}\label{cor:kminus-quotient}
  The quotient norm for $\alg\N/\kminus$ is given by the formula
  \[
    \norm{X+\kminus}=\sup\{\enorm{NXN} \st N\in\N, N<I\}.
  \]
\end{cor}

\begin{proof}

  Clearly for any $K\in\kminus$ and $N<I$ in $\N$,
  $\norm{X-K} \ge \enorm{N(X-K)N} = \enorm{NXN}$,
  and thus $\norm{X+\kminus}$ is at least as big as the supremum of the 
  $\norm{NXN}$'s. We must show the reverse inequality.

  To this end, suppose that $\enorm{NXN}<a$ for all $N<I$, and aim to
  show that $\norm{X+\kminus}\le a$. 
  By \cite{Erdos:OpFiRaNeAl}, $\alg\N$ has a strongly convergent
  approximate identity $E_k$ of finite rank contractions. 
  By \cite[Lemma~4.3]{DavidsonPower:BeApCAl}
  we can assume that $E_k$ converges strong-* and 
  $\lim_n\norm{I-E_n}=1$.
  Choose $N_n$ to be a sequence in $\N$ that increases to $I$,
  and inductively choose compact operators $K_n$ in $\alg\N$ as follows:

  Suppose that $K_1, \ldots, K_{n-1}$ have been chosen to have the property that
  $K_i=K_i(N_i-N_{i-1})$ and
  \[
    \norm{XN_{n-1} - \sum_{i=1}^{n-1} K_i}<a
  \]
  (declaring $N_0 := 0$ for convenience).
  Let $X' := XN_n - \sum_{i=1}^{n-1} K_i$,
  let $P := N_n-N_{n-1}$,
  and observe that the hypotheses of Lemma~\ref{lem:norm-estimate} apply 
  to $X'$, $P$, and $PE_kP$. The result is that
  we can find a $k_0$ for which
  \[
    \norm{XN_n - \sum_{i=1}^{n-1} K_i - XPE_{k_0}P}<a
  \]
  and the induction is completed on taking $K_n := XPE_{k_0}P.$

  Observe that the  series $K := \sum_{i=1}^\infty K_i$ converges weakly, 
  $\norm{X - K}\le a$, and for each $n$, $KN_n = \sum_{i=1}^n K_i$, which is compact,
  hence $K\in\kminus$. Thus $\norm{X+\kminus}\le a$.

\end{proof}

The corresponding result for $\kplus$ follows analogously:

\begin{cor}\label{cor:kplus-quotient}
  The quotient norm for $\alg\N/\kplus$ is given by the formula
  \[
    \norm{X+\kplus}=\sup\{\enorm{N^\perp XN^\perp} \st N\in\N, 0<N\}.
  \]
\end{cor}

We now turn our attention to estimates for the distance from stable ideals not
of compact character. For these, we will need to work with the characterization
of stable ideals in terms of diagonal seminorm functions
(Theorem~\ref{thm:orig-stable-ideal-theorem}).

\begin{lem}\label{lem:dsf}
  Let $\emptydsf$ be a diagonal seminorm function which takes the values
  $\emptydsf=j_N$ for all $0<N<I$. 
  Let $X\in\alg\N$ satisfy $\normn{X}<a$ for all $N$. Then there is a $T\in\alg\N$
  with $\norm{X-T}<a$ and $\normn{T}=0$ for all $N$.
\end{lem}

\begin{proof}

  A routine compactness argument shows that there is a strictly increasing sequence 
  $N_i$ ($i\in\integers$) with
  $\lim_{n\rightarrow-\infty}N_i=0$, $\lim_{n\rightarrow\infty} N_i=I$, and
  \begin{equation}
    \norm{(N_{i+1}-N_{i-1})X(N_{i+1}-N_{i-1})} < a \quad\text{for all $i$.}
    \label{lem:dsf:norm-estimate}
  \end{equation}
  Now consider the possible values of $\emptydsf$ at $N=0$. 
  These are; $0$, $i^+_0$, and $e^+_0$. We shall modify the initial tail of
  the sequence $N_i$ according to which of these values occurs. 
  In the first case, make no changes.
  In the second case (i.e.\ $\dsf{0}{\quad}=i^+_0$) we can renumber the sequence so that
  $N_0 = 0$ and dispense with the negative-index terms,
  while stipulating that estimate~(\ref{lem:dsf:norm-estimate}) still holds.
  In the third case (i.e.\ $\dsf{0}{\quad}=e^+_0$)
  we can again dispense with the negative-index terms,
  as long as we accept that now a weaker estimate holds for the first
  term:
  \begin{equation}\label{lem:dsf:ess-norm-estimate}
    \enorm{(N_2 - N)X(N_2 - N)} < a \quad\text{for all $0<N<N_2$}
  \end{equation}
  (and the original norm estimates still hold for the remaining terms $i\ge 2$).
  However, in this case we can use
  Corollary~\ref{cor:kplus-quotient} to find $K_0=N_2K_0N_2\in\kplus$
  such that $\norm{XN_2-K_0}<a$. The estimate~(\ref{lem:dsf:norm-estimate})
  now holds with $X$ replaced by $X - K_0$ for $i=1$ and for $i\ge 3$. We can recover
  (\ref{lem:dsf:norm-estimate}) for $i=2$ by changing $N_1$ to be a new value
  between $N_0$ and $N_2$ chosen so that $\norm{(N_2-N_1)K_0(N_2-N_1)}$ is
  sufficiently small. Making this change does not affect any of the other norm estimates.

  In the same way we can modify the sequence based on the values of 
  $\emptydsf$ at $N=I$, possibly dispensing with the tail and 
  terminating the sequence at a finite point, $N_{n_0} = I$. 
  In the case that the analogue of
  estimate~(\ref{lem:dsf:ess-norm-estimate}) applies to $N_{n_0}-N_{n_0-2}$,
  we will find a $K_1=N_{n_0-2}^\perp K_0N_{n_0-2}^\perp\in\kminus$
  such that $\norm{N_{n_0-2}^\perp X-K_1}<a$ and adjust $N_{n_0-1}$
  so that estimate~(\ref{lem:dsf:norm-estimate}) still holds for 
  all $i$ for $X - K_1$. 

  Thus, after adjusting our sequence $N_i$ appropriately, 
  estimate~(\ref{lem:dsf:norm-estimate}) will hold for 
  either $X$ itself, or else the operator $X$ adjusted by subtracting
  possibly one or other of $K_0$ and $K_1$. Since~(\ref{lem:dsf:norm-estimate})
  implies that $\norm{N_{i-1}^\perp X N_{i+1}}<a$ for all $i$, it follows
  by a
  version of Arveson's Distance Formula due to Power \cite{Power:DiUpTrOp}, 
  that there is an operator $T$ satisfying
  $N_{i-1}^\perp T N_{i+1}=0$ for all $i$, such that
  $\norm{X-T}<a$.

  Clearly $\normn{T}=j_N(T)=0$ for all $0<N<I$. Next focus on $N=0$. If 
  $\dsf{0}{\quad}=0$ then there is nothing to prove. 
  If $\dsf{0}{\quad}=i^+_0$ then since in this case the sequence of
  $N_i$'s terminates at $N_0=0$ and so $TN_2=0$, therefore $\dsf{0}{T}=0$.
  Lastly, if $\dsf{0}{\quad}=e^+_0$, since $TN_2$ is again $0$ and we
  adjusted $X$ in the previous paragraph by the operator 
  $K_0\in\kplus$ for which $e^+_0(K_0)=0$,
  then we can restore $X$ to its original value and replace $T$ with
  $T+K_0$, observing that $\dsf{0}{T+K_0}$ is also zero.  A similar argument applies 
  for $N=I$ and we are done.
\end{proof}

Recall from \cite[Definition 3.14]{Orr:StIdNeAl} that a diagonal seminorm $\emptydsf$
function is called a greatest diagonal seminorm function if there is
an $X\in\alg\N$ and an $a>0$, so that $\emptydsf$ is the largest diagonal seminorm function
for which $\normn{X}<a$ for all $N\in\N$. Since the collection of all
diagonal seminorm functions is a complete lattice, given any $X$ and $a>0$,
there is a greatest diagonal seminorm function for $X$ and $a$.
It is advantageous to work with
greatest diagonal seminorm functions because of a certain regularity
they exhibit, as shown by the following lemma, quoted from 
\cite[Lemma 3.16]{Orr:StIdNeAl}:

\begin{lem}\label{lem:greatest-diagonal-seminorm-function}
  A diagonal seminorm function, $\emptydsf$, is a greatest diagonal seminorm
  function if and only if it has the following lower semicontinuity property:
  $\emptydsf^-=0$ whenever there is a sequence $N_n$ increasing to $N$ with 
  $\dsf{N_n}{\;\cdot\;}\not=j_{N_n}$ for all $n$, and
  $\emptydsf^+=0$ whenever there is a sequence $N_n$ decreasing to $N$ with 
  $\dsf{N_n}{\;\cdot\;}\not=j_{N_n}$ for all $n$.
\end{lem}

\begin{prp}\label{prp:dsf-estimate}
  Let $\emptydsf$ be a greatest diagonal seminorm function and let
  $X\in\alg\N$ satisfy $\normn{X}<a$ for all $N$. Then there is a $T\in\alg\N$
  with $\norm{X-T}<a$ and $\normn{T}=0$ for all $N$.
\end{prp}

\begin{proof}

  Let $S=\{N \st \emptydsf = j_N\}$. By Lemma~\ref{lem:greatest-diagonal-seminorm-function}, 
  this is an open set and so decomposes as the disjoint union of open intervals $(M_n, N_n)$. 
  For $N\not\in S$ write
  $\emptydsf=a^-_N \vee a^+_N$ where $a^\pm_N$ is one of $0$, $e^\pm_N$, or $i^\pm_N$.
  Note that also by Lemma~\ref{lem:greatest-diagonal-seminorm-function},
  $a^-_{M_n}$ is zero unless it happens that
  $M_n$ is the upper endpoint of another component of $S$ (i.e. $M_n=N_k$ for some $k$).
  Thus it will suffice to construct $T_n=(N_n-M_n)T_n(N_n-M_n)$ with the property
  that $\normn{T_n}=0$ for all $N$ and $\norm{(N_n-M_n)(X-T_n)(N_n-M_n)}<a$, for then
  \[
    T := X - \sum_{n=1}^\infty (N_n-M_n)(X-T_n)(N_n-M_n)
  \]
  will have the desired properties. However if we restrict to the nest
  $(N_n-M_n)\N$ and apply Lemma~\ref{lem:dsf} to $(N_n-M_n)X(N_n-M_n)$ in the
  algebra of this nest, we obtain the desired $T_n$.

\end{proof}

\begin{cor}\label{cor:distance-using-seminorms}
  Let $\I$ be a stable ideal not of compact character
  and let $\F$ be a set of diagonal seminorm functions
  that specify $\I$ as in Theorem~\ref{thm:orig-stable-ideal-theorem}. Then
  the quotient norm is given by the formula
  \[
    \norm{X+\I} = \inf\{\sup_{N\in\N}\normn{X} \st \emptydsf\in\F \}
  \]
\end{cor}

\begin{proof}
  
  If $T\in\I$ and $\e>0$, there is a $\emptydsf$ in $\F$ such that
  $\normn{T}<\e$ for all $N$. Thus
  \[
    \norm{X-T}\ge\normn{X-T} \ge \normn{X}-\e
  \]
  for all $N$. This proves 
  $\norm{X+\I} \ge  \inf\{\sup_{N\in\N}\normn{X} \st \emptydsf\in\F\}$
  and it remains to establish the reverse inequality.

  Suppose that $\inf\{\sup_{N\in\N}\normn{X} \st \emptydsf\in\F\} < a$
  and aim to show that $\norm{X+\I}<a$. It follows there is a
  $\emptydsf\in\F$ such that $\normn{X}<a$ for all $N$. Although
  $\emptydsf$ need not be a \emph{greatest} diagonal seminorm,
  we can take  $\emptydsf'$ to be the greatest diagonal seminorm
  for $X$ and $a$, so that 
  \[
    \normn{X} \le \normn{X}' < a
  \]
  for all $N$. Applying Proposition~\ref{prp:dsf-estimate} 
  to $\emptydsf'$, there is a
  $T\in\alg\N$ with $\norm{X-T}<a$ and 
  $\normn{T}\le\normn{T}'=0$. If follows that $T\in\I$
  and so $\norm{X+\I} <a$.
\end{proof}

\begin{thm}\label{thm:distance-formula}
  Let $\I$ be a stable ideal of $\alg\N$ and let $\O$ be a stable net that
  specifies it in the sense of Theorem~\ref{thm:stable_ideals}. Then
  the quotient norm is given by the formula
  \[
    \norm{X+\I} = \lim_{P\in\O}\sup_{E\in\P} \enorm{EXE}
  \]
\end{thm}

\begin{proof}
  If $\e>0$ we can find a $T\in\I$ such that
  \[
    \norm{X-T} < \norm{X+\I} +\e
  \]
  and so,
  \[
    \lim_{P\in\O}\sup_{E\in P}\enorm{EXE}
      =\lim_{P\in\O}\sup_{E\in P}\enorm{E(X-T)E}
      \le\norm{X-T}
      < \norm{X+\I} +\e
  \]
  Thus $\lim_{P\in\O}\sup_{E\in P}\enorm{EXE}\le\norm{X+\I}$
  and it remains to prove the reverse inequality.

  Let $\O'$ be the set of all collections $P_{T,a}$ 
  as $T$ ranges over $\I$ and $a>0$ (see Remark~\ref{rem:define-induced-stable-net}). 
  In Theorem~\ref{thm:stable_ideals}
  we saw that $\O'$ specifies $\I$. Since $\O$ also specifies $\I$,
  it follows by Theorem~\ref{thm:cofinal-nets} that the two nets $\O$ and $\O'$
  are mutually cofinal. Thus for any $X\in\alg\N$
  \[
    \lim_{P\in\O}\sup_{E\in\P} \enorm{EXE}
      = \lim_{P\in\O'}\sup_{E\in\P} \enorm{EXE}.
  \]

  Let $\F$ be a set of diagonal seminorms that specify $\I$ as
  in Theorem~\ref{thm:orig-stable-ideal-theorem}.
  Suppose that $\lim_{P\in\O'}\sup_{E\in P}\enorm{EXE}<a$, and find a
  $P\in\O'$ such that 
  \[
    \sup_{E\in P}\enorm{EXE}<a.
  \]
  Since $P\in\O'$,
  there is a $T\in\I$ and a $c>0$ such that $P = P_{T, c}$. By rescaling
  $T$, we can assume $c=a$. Next, since $T\in\I$, we can find a diagonal
  seminorm function $\emptydsf$ in $\F$ such that $\normn{T}<a$ for all $N\in\N$.
  Thus by Lemma~\ref{lem:compatible-dsf-and-sets-of-intervals}
  $\emptydsf$ is compatible with $P_{T, a}$ and so also with
  $P_{X, a}$, because this contains $P_{T, a}$. Again, by
  Lemma~\ref{lem:compatible-dsf-and-sets-of-intervals}, we conclude
  that $\normn{X}<a$ for all
  $N$. Thus by Corollary~\ref{cor:distance-using-seminorms}, $\norm{X+\I} \le a$.
  The conclusion is thus that
  \[
    \norm{X+\I} 
    \le \lim_{P\in\O'}\sup_{E\in P}\enorm{EXE}
    = \lim_{P\in\O}\sup_{E\in P}\enorm{EXE}
  \]
  and the result follows.
\end{proof}

\section{Sums of Ideals}\label{sec:sums-of-ideals}

In this section we shall study the algebraic sum of two stable ideals, and show
that this is also stable. In Proposition~\ref{prp:sum-of-stable-ideals}
we present a natural description of the sum
of two stable ideals in terms of the stable nets giving rise to the summands.
  
\begin{lem}\label{lem:sums-are-closed}
  Let $\I_1$ and $\I_2$ be two stable ideals of $\alg\N$. Then the algebraic 
  sum, $\I_1+\I_2$, is a stable ideal.
\end{lem}

\begin{proof}
  Clearly $\I_1+\I_2$ is an automorphism invariant ideal, and we need only show
  that it is norm closed. By a standard argument,
  \[
    \frac{\alg\N}{\I_1}
      \supseteq \frac{\I_1+\I_2}{\I_1}
      \cong \frac{\I_2}{\I_1\cap\I_2}
  \]
  Since $\I_2/\I_1\cap\I_2$ is complete as a normed space, it suffices to show that the
  above algebra isomorphism is isometric, for then $\I_1+\I_2/\I_1$ is closed in 
  $\alg\N/\I_1$, and $\I_1+\I_2$ is its preimage under the quotient map. To this end,
  choose $X\in\I_2$ and let the stable nets $\O_1$ and $\O_2$ respectively induce
  $\I_1$ and $\I_2$. Clearly
  \[
    \norm{X+\I_1} \le \norm{X+(\I_1+\I_2)}
  \]
  On the other hand, by Theorem~\ref{thm:distance-formula}, 
  given $\e>0$, we can find $P_1\in\O_1$ such that
  \[
    \sup_{E\in P_1}\enorm{EXE}<\norm{X+\I_1} + \e
  \]
  and, since $X\in\I_2$, we can certainly find $P_2\in\O_2$ such that
  \[
    \sup_{E\in P_2}\enorm{EXE}<\norm{X+\I_1} + \e
  \]
  But then $P_1\cup P_2\in\O_1+\O_2$ and, by Proposition~\ref{prp:intersect}
  and Theorem~\ref{thm:distance-formula},
  \[
    \begin{split}
      \norm{X+(\I_1\cap\I_2)} 
        & = \lim_{P\in\O_1+\O_2}\sup_{E\in P} \enorm{EXE} \\
        & \le \max_{i=1,2}\sup_{E\in P_i}\enorm{EXE}     \\
        & < \norm{X+\I_1} + \e
    \end{split}
  \]
  Since $\e$ was arbitrary, the result follows.
\end{proof}

\begin{lem}\label{lem:join-of-dsfs}
  Let $\emptydsf^{(1)}$ and  $\emptydsf^{(2)}$ be two diagonal seminorm
  functions and write $\emptydsf^{(1)}\wedge\emptydsf^{(2)}$ for their
  meet in the lattice of diagonal seminorm functions. Suppose that
  $\emptydsf^{(1)}$ is in fact a greatest diagonal seminorm function
  and that for some $Y\in\alg\N$ and $\e>0$
  \[
    \dsf{N}{Y}^{(1)} \wedge \dsf{N}{Y}^{(2)} < \e \text{ for all } N\in\N
  \]
  Then $Y=Y_1+Y_2$ for some $Y_i\in\alg\N$ where
  \[
    \dsf{N}{Y_1}^{(1)} < \e \text{ and } \dsf{N}{Y_2}^{(2)} < \e \text{ for all } N\in\N
  \]
\end{lem}

\begin{proof}
  By Lemma~\ref{lem:greatest-diagonal-seminorm-function},
  $S := \{N\st \emptydsf^{(1)} = j_N\}$
  is an open set, which decomposes as the disjoint union of intervals
  $(M_n, N_n)$. Let $Y_1 := Y - \sum_n (N_n-M_n)Y(N_n-M_n)$. Clearly
  $\dsf{N}{Y_1}^{(1)}=0$ for all $N\in S$ and observe that
  for $N\not\in S$, 
  $\emptydsf=a^-_N \vee a^+_N$ where $a^\pm_N$ is one of $0$, $e^\pm_N$, or $i^\pm_N$.
  Also by Lemma~\ref{lem:greatest-diagonal-seminorm-function},
  $a^-_{M_n}$ (respectively, $a^+_{N_n}$) is zero unless it happens that
  $M_n$ is the upper endpoint (respectively, $M_n$ is the lower endpoint)
  of another component of $S$, and so $\dsf{N}{Y_1}^{(1)}=0$ for all $N\in\N$. 

  Next, consider $Y' := Y - Y_1 = \sum_n (N_n-M_n)Y(N_n-M_n)$. 
  If $N\in S$ then $\emptydsf^{(1)} = j_N \ge \emptydsf^{(2)}$ and so,
  since $\dsf{N}{Y}^{(i)} \ge \dsf{N}{Y'}^{(i)}$, we cannot have 
  $\dsf{N}{Y'}^{(2)}\ge \e$. Thus $\dsf{N}{Y'}^{(2)} < \e$ for all $N\in S$.

  For $N\not\in S$ write $\emptydsf^{(i)-}$ and $\emptydsf^{(i)+}$ for 
  the left and right respective parts of $\emptydsf^{(i)}$. In other words, 
  $\emptydsf^{(i)} = \emptydsf^{(i)-}\vee\emptydsf^{(i)+}$ and each 
  $\emptydsf^{(i)\pm}$ is one of $0$, $e_N^\pm$, $i_N^\pm$.
  Let $A^- := \{n \st \dsf{N_n}{Y'}^{(2)-} \ge \e \}$ and
  $A^+ := \{n \st \dsf{M_n}{Y'}^{(2)+} \ge \e \}$.
  If $n\in A^-$ then 
  \[
    \dsf{N_n}{Y'}^{(1)-} < \e \le \dsf{N_n}{Y'}^{(2)-}
  \]
  Conversely, any $N\in(M_n, N_n)$ that satisfies $\dsf{N}{Y'}^{(1)}\ge \e$,
  must have $\dsf{N}{Y'}^{(2)}< \e$, and so $\emptydsf^{(2)}\not= j_N$
  (since $j_N$ dominates all other diagonal seminorms).
  It follows by Lemma~\ref{lem:greatest-diagonal-seminorm-function}
  that since $\dsf{N_n}{Y'}^{(2)-}\not=0$, $N_n$ cannot be
  the limit from below of such $N$. Thus there is an $L_n\in(M_n, N_n)$
  such that $\dsf{N}{Y'}^{(1)}< \e$ for all $N\in[L_n, N_n)$.

  By the same token, if $n\in A^+$
  we can find $G_n\in(M_n, N_n)$ so that $\dsf{N}{Y'}^{(1)}< \e$ 
  for all $N\in(M_n, G_n]$. Without loss, insist that
  $M_n<G_n<L_n<N_n$ for $n\in A^-\cap A^+$.

  Let 
  $Y'_1 := \sum_{i\in A^-} (N_i - L_i)Y(N_i - L_i) + \sum_{j\in A^+} (G_j - M_j)Y(G_j - M_j)$.
  Then $\dsf{N}{Y_1'}^{(1)}<\e$ 
  for all $N$, by a similar argument to the first paragraph,
  and clearly  $\dsf{N}{Y' - Y_1'}^{(2)}<\e$ for all $N$.
  Thus $Y=(Y_1 + Y_1') + (Y' - Y_1')$, which has the desired properties.
\end{proof}

\begin{prp}\label{prp:sum-of-stable-ideals}
  Let $\I_1$ and $\I_2$ be stable ideals specified respectively by the
  stable nets $\O_1$ and $\O_2$. For $P_1\in\O_1$ and $P_2\in\O_2$ define
  \[
    P_1\cdot P_2 := \{E_1 E_2 \st E_1\in P_1, E_2\in P_2\}
  \]
  and then define
  \[
    \O_1\cdot\O_2 := \{P_1\cdot P_2 \st P_1\in\O_1, P_2\in O_2\}
  \]
  Then $\O := \O_1\cdot\O_2$ is a stable net, and
  \[
    \I_1 + \I_2
    =
    \{X\in\alg\N \st \lim_{P\in\O}\sup_{E\in P}\enorm{EXE} = 0\}
  \]
\end{prp}

\begin{proof}
  It is routine to verify that $\O:=\O_1\cdot\O_2$ is a stable net
  and whenever $Y=Y_1+Y_2\in\I_1+\I_2$ then
  $\lim_{P\in\O}\sup_{E\in P}\enorm{EYE}=0$. It remains to suppose that 
  $\lim_{P\in\O}\sup_{E\in P}\enorm{EYE}=0$ for some $Y\in\alg\N$
  and show that $Y\in\I_1+\I_2$.

  Since the ideals $\I_i$ ($i=1,2$) are stable ideals, there are families of
  diagonal seminorm functions $\F_i$ which determine them. Let $\O'_i$
  be the set of $P_{T,a}$ as $T$ varies over $\I_i$ and $a>0$. As we saw in
  Theorem~\ref{thm:stable_ideals}, each $\O'_i$ determines $\I_i$ and so,
  by Theorem~\ref{thm:cofinal-nets}, each $\O'_i$ is mutually cofinal
  with each $\O_i$. Thus it is easily seen that $\O':=\O'_1\cdot\O'_2$
  and $\O$ are mutually cofinal.

  Given $\e>0$, since 
  $\lim_{P\in\O'}\sup_{E\in P}\enorm{EYE} = \lim_{P\in\O}\sup_{E\in P}\enorm{EYE}=0$, 
  we can find a $P\in\O'$ such that $\enorm{EYE}<\e$ for all $E\in P$.
  However $P=P_1\cdot P_2$ where each $P_i = P_{T_i,\e}$ for a suitably
  scaled $T_i\in\I_i$. Thus we can find $\emptydsf^{(i)}\in\F_i$ such that
  $\dsf{N}{T_i}^{(i)}<\e$ for all $N$. Since the proof of 
  Lemma~\ref{lem:join-of-dsfs} requires $\emptydsf^{(1)}$ to be a greatest
  diagonal seminorm function, let us replace $\emptydsf^{(1)}$ with
  the greatest diagonal seminorm function for which 
  $\dsf{N}{T_1}^{(1)}<\e$ for all $N$.

  By definition, each $P_i$ is compatible with $\emptydsf^{(i)}$
  in the sense of Definition~\ref{def:compatible-collections}.
  It is straightforward to check that therefore 
  $P=P_1\cdot P_2$ is compatible with $\emptydsf^{(1)}\wedge\emptydsf^{(2)}$.
  Therefore also $P_{Y, \e}\supseteq P$ is compatible with 
  $\emptydsf^{(1)}\wedge\emptydsf^{(2)}$and so,
  by Lemma~\ref{lem:compatible-dsf-and-sets-of-intervals}, 
  $\dsf{N}{Y}^{(1)}\wedge\dsf{N}{Y}^{(2)}<\e$.

  It follows from Lemma~\ref{lem:join-of-dsfs} that $Y=Y_1+Y_2$
  where $\dsf{N}{Y_i}^{(i)} < \e$ for all $N$ and $i=1,2$. 
  Although for technical reasons we replaced the original
  $\emptydsf^{(1)}$ from $\F_1$ with the greatest diagonal seminorm 
  function for $T_1$ and $\e$, nevertheless since that seminorm function
  dominated the original one, the current estimate for $Y_1$
  holds with the original $\emptydsf^{(1)}$ from $\F_1$ too.
  But thus, by Theorem~\ref{thm:distance-formula}, $Y$ is at
  most distance $2\e$ from $\I_1+\I_2$. Finally, since
  $\e$ was arbitrary, Lemma~\ref{lem:sums-are-closed}
  shows that $Y\in\I_1+\I_2$.
\end{proof}

\bibliographystyle{abbrv}
\bibliography{bibliography}

\end{document}